\documentclass[%
11pt,
reprint,
onecolumn,
tightenlines,
superscriptaddress,
preprintnumbers,
nofootinbib,
amsmath,amssymb,amsthm,
physrev,
eqsecnum,tikz,
]{revtex4-2}

\usepackage{isomath}
\usepackage{amsmath,amsthm}
\usepackage{amsbsy}
\usepackage{amssymb}
\usepackage{amscd}
\usepackage{amsfonts}
\usepackage{stmaryrd}
\usepackage{siunitx}
\usepackage{euscript}
\usepackage[utf8]{inputenc}
\usepackage[T1]{fontenc}
\usepackage{newtxtext} 
\everymath{\displaystyle}
\usepackage{exscale}
\usepackage{microtype}
\usepackage{hyperref}
\usepackage{booktabs}
\usepackage{algorithm}
\usepackage{algpseudocode}

\usepackage{graphicx}
\usepackage{boxedminipage}
\usepackage{calc}
\usepackage[usenames,dvipsnames]{xcolor}
\graphicspath{ {media/} }
\usepackage[caption=false,justification=centerlast]{subfig}

\usepackage{setspace}
\usepackage{enumitem}
\setitemize{noitemsep,topsep=0pt,parsep=0pt,partopsep=0pt}
\setenumerate{noitemsep,topsep=0pt,parsep=0pt,partopsep=0pt}
\setdescription{noitemsep,topsep=0pt,parsep=0pt,partopsep=0pt}

\usepackage[small]{titlesec}

\titlespacing*{\section}{0pt}{12pt plus 4pt minus 2pt}{2pt plus 2pt minus 2pt}
\titlespacing*{\subsection}{0pt}{12pt plus 4pt minus 2pt}{2pt plus 2pt minus 2pt}
\titlespacing*\subsubsection{0pt}{12pt plus 4pt minus 2pt}{2pt plus 2pt minus 2pt}
\titlespacing*\paragraph{0pt}{12pt plus 4pt minus 2pt}{2pt plus 2pt minus 2pt}

\makeatletter

    \renewcommand*{\p@subsection}{}
    
    \renewcommand*{\p@subsubsection}{}
\makeatother

\usepackage{isomath}
\usepackage{amsmath}
\usepackage{amssymb}
\usepackage{amscd}
\usepackage{amsfonts}
\usepackage{upgreek}

\theoremstyle{definition}

\newtheorem{remark}{Remark}[section]
\AtEndEnvironment{remark}{\null\hfill\qedsymbol}
\newtheorem{example}{Example}[section]
\AtEndEnvironment{example}{\null\hfill\qedsymbol}

\newcommand{\bfchi}{\mathbold {\chi}}

\newcommand{\bfsigma}{\mathbold {\sigma}}

\newcommand{\bfzero}{\mathbf{0}}

\DeclareMathOperator{\divergence}{div}

\DeclareMathOperator{\trace}{tr}

\newcommand{\parderiv}[2]{\frac{\partial #1}{\partial #2}}
\newcommand{\dm}{\ \mathrm{d}}
\newcommand{\deriv}[2]{\frac{\dm #1}{\dm #2}}

\newcommand{\derivM}[2]{\frac{\ \mathrm{D} #1}{\ \mathrm{D} #2}}

\newcommand{\bfb}{{\mathbold b}}

\newcommand{\bfe}{{\mathbold e}}

\newcommand{\bfg}{{\mathbold g}}

\newcommand{\bfm}{{\mathbold m}}
\newcommand{\bfn}{{\mathbold n}}

\newcommand{\bfq}{{\mathbold q}}

\newcommand{\bft}{{\mathbold t}}
\newcommand{\bfu}{{\mathbold u}}
\newcommand{\bfv}{{\mathbold v}}

\newcommand{\bfx}{{\mathbold x}}

\newcommand{\bfF}{{\mathbold F}}

\newcommand{\bfI}{{\mathbold I}}

\newcommand{\bfP}{{\mathbold P}}

\newcommand{\bfV}{{\mathbold V}}

\newcommand{\bfX}{{\mathbold X}}

\begin{document}

\preprint{To appear in Journal of the Mechanics and Physics of Solids (DOI: \href{https://doi.org/10.1016/j.jmps.2025.106076}{10.1016/j.jmps.2025.106076})
.}

\title{Accretion and Ablation in Deformable Solids using an Eulerian Formulation:\\ A Finite Deformation Numerical Method}

\author{S. Kiana Naghibzadeh}
    \email{kiana@mit.edu}
    \affiliation{Department of Civil and Environmental Engineering, Massachusetts Institute of Technology}

\author{Anthony Rollett}
    \affiliation{Department of Materials Science and Engineering, Carnegie Mellon University}

\author{Noel Walkington}
    \affiliation{Center for Nonlinear Analysis, Department of Mathematical Sciences, Carnegie Mellon University}

\author{Kaushik Dayal}
    \affiliation{Department of Civil and Environmental Engineering, Carnegie Mellon University}
    \affiliation{Center for Nonlinear Analysis, Department of Mathematical Sciences, Carnegie Mellon University}
    \affiliation{Department of Mechanical Engineering, Carnegie Mellon University}

\date{\today}

\begin{abstract}
    Surface growth, i.e., the addition or removal of mass from the boundary of a solid body, occurs in a wide range of processes, including the growth of biological tissues, solidification and melting, and additive manufacturing. 
    To understand nonlinear phenomena such as failure and morphological instabilities in these systems, accurate numerical models are required to study the interaction between mass addition and stress in complex geometrical and physical settings.
    Despite recent progress in the formulation of models of surface growth of deformable solids, current numerical approaches require several simplifying assumptions.
    
    This work formulates a method that couples an Eulerian surface growth description to a phase-field approach.
    It further develops a finite element implementation to solve the model numerically using a fixed computational domain with a fixed discretization.
    This approach bypasses the challenges that arise in a Lagrangian approach, such as having to construct a four-dimensional reference configuration, remeshing, and/or changing the computational domain over the course of the numerical solution.
    It also enables the modeling of several settings --- such as non-normal growth of biological tissues and stress-induced growth --- which can be challenging for available methods.
    
    The numerical approach is demonstrated on a model problem that shows non-normal growth, wherein growth occurs by the motion of the surface in a direction that is not parallel to the normal of the surface, that can occur in hard biological tissues such as nails, horns, etc.
    Next, a thermomechanical model is formulated and used to investigate the kinetics of freezing and melting in ice under complex stress states, particularly to capture regelation which is a key process in frost heave and basal slip in glaciers.
\end{abstract}

\maketitle
\section{Introduction}

Surface growth is the addition (accretion) or removal (ablation) of mass from the boundary of a solid body.
Numerous natural and engineering processes involve interactions between surface growth and stress, and this interaction plays a major role in both the mechanisms of surface growth as well as the development of stresses within the growing deformable solid body. 
Notable examples include the melting and refreezing of glacial ice, where thermomechanical stresses give rise to crack (crevasse) propagation \cite{colgan2016glacier} and basal pressure melting and regelation \cite{cuffey2010physics}; 
additive manufacturing, where thermomechanical stresses introduce residual stress and distortion in printed objects \cite{fitzharris2018effects}; 
electrochemical deposition in batteries, where instabilities and dendrite growth can be significantly influenced by external pressure \cite{yin2018insights}; 
stress-induced growth of biological tissues \cite{ambrosi2007stress}; and cell motility \cite{salvadori2024generation}. 

A continuum description of surface growth requires a time-evolving set of material particles, in which both the current and reference configurations of the growing body are  evolving in time. 
This feature presents significant challenges both in terms of formulating the problem within the framework of continuum mechanics and in devising numerical methods to solve it. 

\paragraph*{Prior Work on Formulation of Surface Growth.} 

Solid mechanics typically uses the Lagrangian description, i.e., material particles are labeled by, and the deformation is computed using, the locations of the material particles in the reference configuration.
This description is convenient when the set of material particles remains constant, which is not the case when dealing with surface growth phenomena.
The extension of Lagrangian approaches to problems of surface growth typically include a time-evolving reference configuration, to incorporate the history of particle attachment to the object \cite{ateshian2007theory, hodge2010continuum, tomassetti2016steady, sozio2017nonlinear, abi2019kinetics, pradhan2023accretion}. 
This introduces significant complexity but does allow for the use of standard methods to compute the deformation gradient $\bfF$ that is required for the stress response. 
In the limit of infinitesimal deformations, the reference configuration can be approximated with the current configuration, making the formulation of surface growth less challenging \cite{naumov1994mechanics, kadish2008stresses}, but this is obviously limited in applicability. 

An alternate approach to surface growth is the use of Eulerian-based approaches \cite{naghibzadeh2021surface, ganghoffer2018combined}, that work directly with the current configuration and bypass the direct definition of the reference configuration, e.g. \cite{clayton2013nonlinear, clayton2019nonlinear}. 
However, this raises the challenge of computing the deformation gradient $\bfF$ --- that is required for the stress response --- when we do not have an explicit description of the reference configuration.
In our previous work, our approach to this problem was composed of 2 key elements \cite{naghibzadeh2021surface,naghibzadeh2022accretion}.
First, we use an evolution equation for $\bfF$ that is posed in the current configuration, following \cite{liu2001eulerian}.
Second, we introduced a new kinematic descriptor, the \textit{elastic deformation}, $\bfF_e$, that accounts for the part of $\bfF$ that causes stress.
Using these elements together, we can completely eliminate $\bfF$ from the formulation in favor of $\bfF_e$, where the time-evolution of the latter is governed by a transport-type PDE.
The key advantages of using $\bfF_e$ are: first, it is not constrained to satisfy the compatibility constraint, unlike $\bfF$; and, second, for the material particles that are being attached to the body during surface growth, the reference configuration --- and consequently $\bfF$ --- are not well-defined quantities, whereas the stress --- and consequently $\bfF_e$ --- are well-defined.

All of these elements provide an approach wherein we can compute the stress and deformation of a growing body without explicitly computing or using the complex reference configuration.
Unlike previously-developed Lagrangian approaches, this approach is not limited by model-specific assumptions required for extending the reference configuration; a detailed comparison of this method with some of the available Lagrangian methods is discussed in \cite{naghibzadeh2022accretion}.

\paragraph*{Prior Work on Numerical Methods.} 

Developing the formulation and numerical modeling of the thermomechanics of surface growth continues to be an active area of research, e.g. \cite{pradhan2024nonlinear, li2024mechanical}. 
Recently, \cite{von2021morphogenesis} developed a Lagrangian finite element (FE) method to study surface growth based on the chemo-mechanical growth formulation developed in \cite{abi2019kinetics} for simple geometries with uniform growth velocity.
    A challenge in this approach is that the Lagrangian formulation requires potentially-expensive extension and remeshing of the computational domain at each timestep. 
    Further, in fluid-solid interaction problems as well as situations where solids experience very large deformation --- e.g., at high rates --- standard Lagrangian methods are often fail because of the extreme mesh distortions \cite{zukas2004introduction,anderson1987overview,benson1992computational}.

To address these challenges, several Eulerian \cite{van1994eulerian, dunne2006eulerian, kamrin2012reference, kamrin2009eulerian, sugiyama2011full, miller2001high, liu2001eulerian, okazawa2007eulerian} and hybrid Eulerian/Lagrangian \cite{hu2001direct, de2020material, belytschko2014nonlinear, madadi2024subdivision, moutsanidis2019modeling} numerical methods have emerged.
    However, we highlight that these methods consider problems in which the set of material particles in the solid is preserved in time, and hence they are based on convecting material points or using the mapping between the current and reference configuration.
Consequently, they are not applicable to modeling surface growth problems involving evolving sets of material particles. 
Recently, \cite{bergel2021finite} extended a combination of the Arbitrary Lagrangian-Eulerian FE method and the updated Lagrangian method \cite{belytschko2014nonlinear} to study surface growth for discrete accretion of material using a moving mesh, but are restricted to the case that the added material is stress-free.

\paragraph*{Contributions of This Paper.} 

In this paper, we focus on developing a fully-Eulerian FE framework based on the formulation developed in \cite{naghibzadeh2021surface} to numerically simulate surface growth for arbitrary growth velocity, geometry, and stress state of the added material (but neglecting the effect of inertia). 
We exploit the advantage of an Eulerian formulation to solve the governing equations on a fixed computational domain with a fixed discretization. 
This requires the use of a domain larger than the initial body to contain the growing deformable solid body in the computational domain for the entire time interval of interest. 
Therefore, we use a phase function to distinguish the solid body from the rest of the computational domain, following ideas that have been used extensively in the numerical simulation of multiphase problems, e.g. \cite{osher2005level, tryggvason2007immersed,hou1999level}.

We apply the numerical method to study 2 examples.
First, we demonstrate model simulations of ``non-normal growth'' \cite{skalak1982analytical,skalak1997kinematics}, observed in several hard biological tissues such as horn and nail.
Here, growth occurs at a fixed surface and the generating cells that push the growing material outward are not aligned with the normal vector of the surface \cite{skalak1997kinematics}. 
Hence, the growth of a stress-free material in a non-normal direction is observed.
However, conventional interface-based methods use only the normal interface velocity and cannot directly simulate non-normal growth; we use the elastic deformation to model non-normal growth by applying a pre-stress to the accreted material that drives shear as the pre-stress relaxes.

Second, we develop a thermomechanical extension of the phase-field model to study stress-induced melting and refreezing (i.e., regelation) in ice, based on ideas from \cite{abeyaratne2006evolution,penrose1990thermodynamically,hou1999level,fried1994dynamic,agrawal2015dynamic,agrawal2015dynamic-2,ateshian2022continuum}.
The typical approach to studying this phenomenon is by using the Clausius-Clapeyron relation, which relates the melting temperature to the hydrostatic pressure in the solid.
However, the Clausius-Clapeyron relation is difficult to incorporate in a consistent thermomechanical setting, wherein the stress is obtained from solving the field equation of momentum balance and is hence typically not hydrostatic \cite{sekerka2004solid, style2023generalized}. 
The free energy approach does not require the explicit use of the Clausius-Clapeyron relation, but it is instead an outcome of the formulation of the free energy, the balance laws, and a kinetic law.
This approach is consequently applicable to realistic problems with spatially-varying non-hydrostatic stress fields.
The problem of regelation also involves both fluid and solid phases with transformations between the phases in both directions; the deformation of the fluid phase would, in general settings, be exceptionally challenging for a Lagrangian formulation.

\paragraph*{Organization.}

Section \ref{sec:summaryOfFormularion} summarizes the Eulerian formulation of the surface growth problem developed in \cite{naghibzadeh2021surface}. 
Section \ref{sec:NumericalAlorithm} develops the numerical method to solve the quasi-static growth formulation. 
Section \ref{sec:examples} uses the method to demonstrate non-normal growth, motivated by observations in some hard biological tissues.
Section \ref{sec:MeltingUnderPressure} formulates a thermomechanical model of the ice-water phase transformation, and demonstrates stress-driven melting and refreezing using the example of the passage of a wire through a block of ice without splitting it \cite{bottomley1872melting}.

\section{Formulation of the Growth Problem}
\label{sec:summaryOfFormularion}

\begin{figure}[htb!]
    \includegraphics[width=0.6\textwidth]{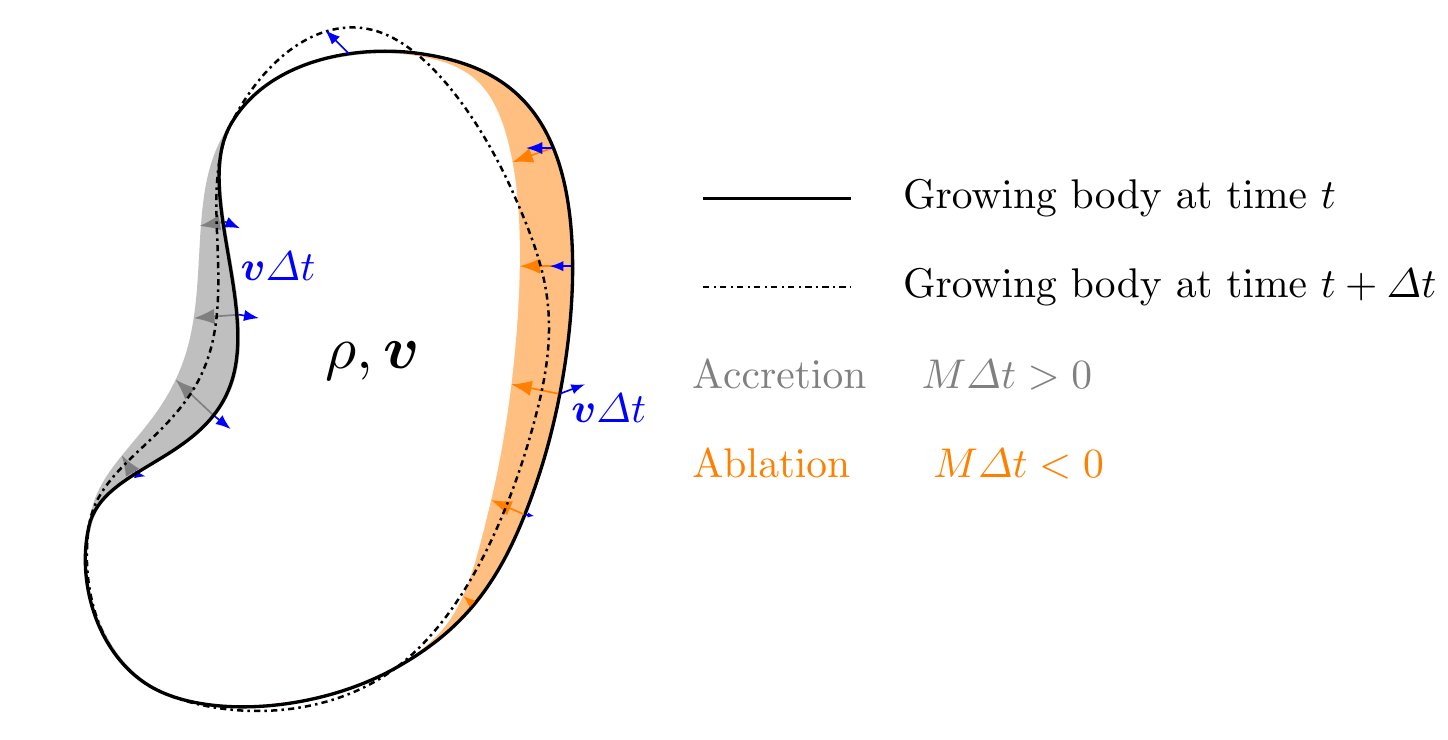}
    \caption{A schematic of the evolution of a growing body in a small interval of time $\Delta t$. The spatial location of the growing boundary depends on the growth velocity and the continuum particle velocity at the boundary.}
    \label{fig:definitions}
\end{figure}

We model the process of surface growth through the introduction of sources of mass and linear momentum that are localized on the growing surface, following \cite{naghibzadeh2021surface}.
To define the source terms, we specify the rate of addition of mass per unit area $M$ and linear momentum per unit area $\bfP$.
Accretion is modeled by $M>0$ and ablation by $M<0$ (Fig. \ref{fig:definitions}).
From these quantities, we can infer the velocity (or specific momentum), $\bfv_a := \bfP / M$, of the added material at the instant of attachment.
We assume that the added particles do not carry angular momentum, such as due to individual particle spins; therefore, the balance of angular momentum provides only the standard result that the Cauchy stress must be symmetric.

\paragraph*{Notation.}
All integration and differentiation operations, e.g., divergence and gradient operators, are with respect to spatial coordinates unless explicitly noted.
$\bfsigma$ is the Cauchy stress; $\rho$ is the mass density; $\bfv$ is the particle velocity; $\bfb$ is the body force; and $\bfV_b$ is the velocity of the growing surface.

\subsection{Balance Laws}
\label{sec:balances}

\begin{figure}[htb!]
    \centering
    \includegraphics[width=0.9\textwidth]{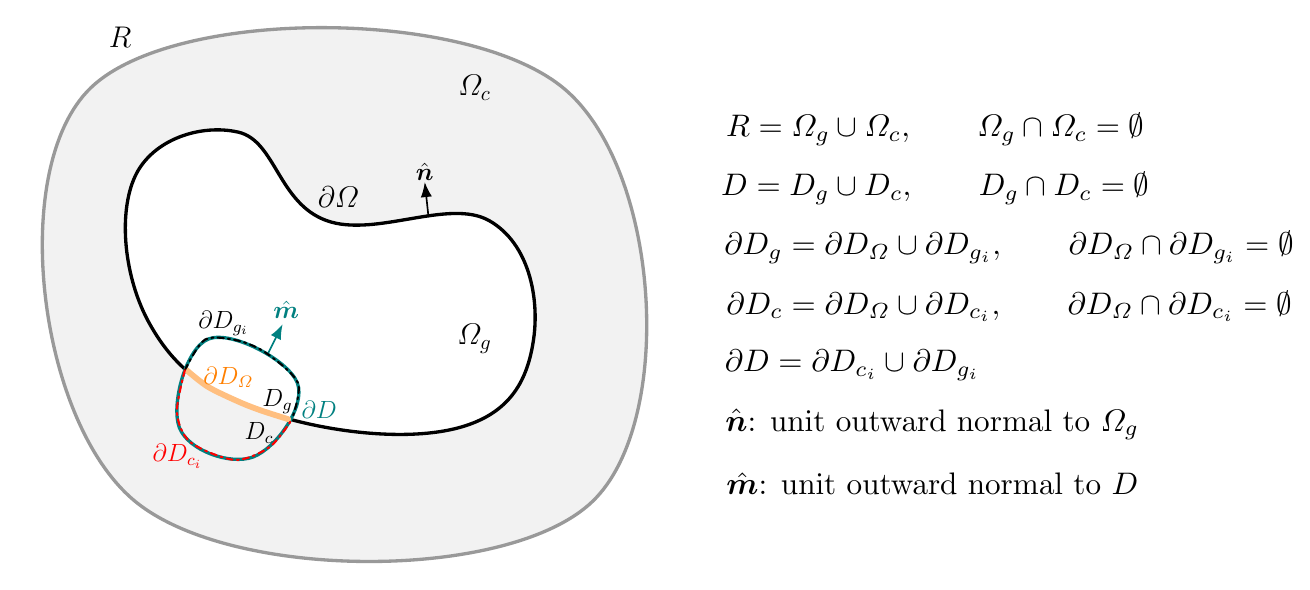}
    \caption{Schematic of a growing body $\Omega_g$ and the ambient exterior $\Omega_c $. The balance equations are derived using the arbitrary subdomain $D$.}
    \label{fig:gov-eq}
\end{figure}

We use the Eulerian description to avoid explicitly introducing a time-dependent reference configuration with a time-dependent set of material points.
The balance laws are completely standard, and the main difference is the presence of source terms in the jump conditions at the growing surface.
The source terms can be understood as a coarse-grained approach to treating the complex process of growth without considering any of the details but only specifying it in terms of net mass, momentum, and energy transfer at the boundary.
Further, when the ambient environment outside of the body at the growing surface has a negligible effect on the mechanics of the body, the domain now involves only the solid body and the jump conditions become boundary conditions (Section \ref{sec:neglect-exterior}).

Figure \ref{fig:gov-eq} is a schematic that defines our notation.
The domain $R(t)$ consists of the growing body $\Omega_g(t)$ with boundary $\partial \Omega (t)$, and the ambient exterior $\Omega_c (t) = R(t) \setminus \Omega_g(t)$. 
An arbitrary subdomain $D(t)$ inside $R(t)$ will be used to find the local forms of the balance laws; when $D(t)$ intersects $\partial \Omega (t)$, we will obtain the jump conditions.
The singular source terms due to growth will be localized on  $\partial\Omega$.
We further choose $D$ such that its boundary velocity is equal to the particle velocity at the boundary; therefore, there is no flux at the boundary of $D$ except due to growth.

We use $\llbracket \alpha \rrbracket := \alpha^g - \alpha^c$ to denote the jump operator, where $\alpha^g$ and $\alpha^c$ are the limits of the discontinuous variable $\alpha$ when approaching the surface of discontinuity from within $D_g$ and within $D_c$, respectively.

\subsubsection{Mass Balance}

When $D$ does not intersect $\partial \Omega $, the balance of mass has the form:
\begin{equation*}
    \deriv{\ }{t} \int_D \rho \dm V = 0
    \implies
    \parderiv{\rho}{t} + \divergence(\rho \bfv) = 0 
\end{equation*}

However, if $D$ intersects $\partial \Omega$, the total mass inside $D$ is not constant. 
The rate of addition of mass is $M$, and therefore, the mass balance equation is:
 \begin{equation}
 \label{eqn:mass-balace-bulk}
    \deriv{\ }{t} \int_{D} \rho \dm V =
    \int_{\partial D_{\Omega} } M \dm A
\end{equation}

Following, e.g., \cite{abeyaratne1952lecture}, gives the jump condition at the singular surface:
\begin{equation}
    \llbracket \rho (\bfV_b - \bfv) \cdot \hat{\bfn} \rrbracket  = M
    \quad \text{ on } \partial \Omega
\end{equation}

\subsubsection{Momentum Balance}

When $D$ does not intersect $\partial \Omega$, the balance of momentum has the form:
\begin{equation}
\label{eqn:mom-bal-bulk-sum}
    \deriv{\ }{t} \int_{D} \rho \bfv \dm V
    =
    \int_{D} \rho \bfb \dm V + \int_{\partial D} \bft \dm A 
    \implies
    \parderiv{\ }{t} (\rho \bfv) +
    \divergence (\rho \bfv \otimes \bfv)=
    \rho \bfb +
    \divergence(\bfsigma)
\end{equation}
using $\bft = \bfsigma \hat\bfm$, where $\bft$ is the traction vector, $\bfsigma$ is the Cauchy stress tensor, and $\hat\bfm$ is the unit outward normal to $\partial D$ (Figure \ref{fig:gov-eq}).

If $D$ intersects $\partial \Omega$, the material added due to growth can also transport momentum, given by $\hat{P} = M \hat{\bfv}_a$ per unit area of the singular growth surface.
The balance law is then given by:
\begin{equation}
\label{eqn:lin-mom-int}
    \deriv{\ }{t} \int_{D} \rho \bfv \dm V = 
    \int_{D} \rho \bfb \dm V +
    \int_{\partial D} \bft \dm A +
    \int_{\partial D_{\Omega}} M \hat{\bfv}_a \dm A 
\end{equation}
The corresponding jump condition is:
\begin{equation}
    \llbracket \rho \bfv ( (\bfV_b - \bfv) \cdot \hat{\bfn} )\rrbracket + 
    \llbracket \bfsigma \hat{\bfn} \rrbracket = 
    M \hat{\bfv}_a
    \quad \text{ at } \partial \Omega
\end{equation}

Applying the balance of angular momentum when $D$ does not intersect with interface provides the standard result that $\bfsigma$ must be symmetric.
When $D$ includes the growing surface,  the contribution to the angular momentum from the added material is $\bfx \times (M \hat{\bfv}_a)$; however, this is simply ($\bfx \times \eqref{eqn:lin-mom-int})$.
Hence, the balance of linear momentum implies the balance of angular momentum when the sources are considered.
We note that we have ignored the situation when the added material carry intrinsic angular momentum, e.g., due to individual particle spins.

\subsubsection{Neglecting the Exterior Ambient Environment}
\label{sec:neglect-exterior}

When the ambient environment outside of the growing body has a negligible effect on the mechanics of the body, the problem can be significantly simplified.
In particular, the domain of the problem now involves only the growing solid body, and the jump conditions become boundary conditions on $\partial\Omega$:
\begin{align}
    \label{eqn:mass-bc-main}
    \text{Mass: } & \rho (\bfV_b - \bfv) \cdot {\hat \bfn}  = M
    \quad \iff \quad
    \bfV_b \cdot {\hat \bfn} = \bfv \cdot {\hat \bfn}  + \frac{M}{\rho}
    \\
    \label{eqn:mom-bc-main}
    \text{Momentum: } & \rho \bfv (\bfV_b - \bfv) \cdot {\hat \bfn} + 
    \bfsigma \hat {\bfn} - \bft_b = 
    M \bfv_a 
    \quad \iff \quad
    \bfsigma {\hat \bfn} 
    =
    M (\bfv_a - \bfv) + \bft_b
\end{align}
where $\bft_b$ is the traction due to the external forces at the boundary of the growing body, and we have used $M = \rho (\bfV_b - \bfv) \cdot \bfn$.

Following \cite{skalak1997kinematics}, we define the growth velocity $\bfv_g$ at the boundary such that $\bfv_g \cdot \hat \bfn = \frac{M}{\rho}$, and hence \eqref{eqn:mass-bc-main} can be written as $\bfV_b \cdot {\hat \bfn} = \bfv \cdot {\hat \bfn}  + \bfv_g \cdot {\hat \bfn}$.
We notice here that only the normal components of $\bfV_b$ and $\bfv_g$ are of physical significance.

In the case of \textit{slow growth}, the inertial terms $\rho \bfv (\bfV_b - \bfv) \cdot \hat{\bfn}$ and $M \hat{\bfv}_a$ are negligible relative to $\bfsigma$, and the traction boundary condition reduces to $\bfsigma \hat{\bfn} = \bft_b$.

\subsection{Kinematics: Evolution of the Elastic Deformation}
\label{sec:Kinematic}

The stress response of solid materials generally requires knowledge of the deformation gradient $\bfF = \parderiv{\bfchi}{\bfX}( \bfX , t)$, 
where $\bfX$ is the location in the reference configuration of a material point that is undergoing the motion $\bfx = \bfchi(\bfX , t)$. 
To find $\bfF$, we adapt the approach from \cite{liu2001eulerian} that uses the following first-order linear transport/advection equation that governs the time-evolution of $\bfF$.
We can write the material time derivative of $\bfF$ as:
\begin{align}
    & \parderiv{\bfF(\bfX,t)}{t}
    =
    \parderiv{\ }{t}\parderiv{\bfchi(\bfX,t)}{\bfX}
    =
    \parderiv{\ }{\bfX}\parderiv{\bfx}{t}
    =
    \parderiv{\bfv}{\bfX}
    =
    \parderiv{\bfv}{\bfx}\parderiv{\bfx}{\bfX}
    =
    \left(\nabla\bfv\right) \bfF
    \\
    &
    \Rightarrow
    \dot{\bfF}=\nabla\bfv \bfF
\end{align}
where $\dot{\bfg}$ represents the material time derivative of a quantity $\bfg$.
Writing this in Eulerian form, we have:
\begin{equation}
    \label{eqn:def-grad-transport}
    \parderiv{\bfF(\bfx,t)}{t}+(\bfv\cdot\nabla)\bfF=\left(\nabla\bfv \right) \bfF
\end{equation}
For well-posedness of \eqref{eqn:def-grad-transport},  we require that boundary conditions be imposed at inlets, where the transport velocity points into the domain, but not at outlets, where the transport velocity points out of the domain \cite{trangenstein2009numerical,naghibzadeh2022accretion}; these correspond to accretion and ablation, respectively.
This introduces a key difficulty in directly using \eqref{eqn:def-grad-transport}: $\bfF$ must satisfy kinematic compatibility, and it is not clear to us how to prescribe the boundary conditions to satisfy this constraint at all times.
Further, while \eqref{eqn:def-grad-transport} preserves the compatibility of a field in time as it evolves, it is not clear to us how to construct a numerical approximation method that is guaranteed to similarly preserve compatibility.

To avoid the difficulties associated with satisfying kinematic compatibility, we begin by noticing that the stress state of the added material, denoted by $\bfsigma^*$, is a well-defined and physically-meaningful quantity, whereas the deformation gradient $\bfF$ corresponding to $\bfsigma^*$ can be chosen for convenience by appropriately defining the reference at the time of attachment.
Hence, we introduce two kinematic fields: (1) $\bfF_{relax}$ to quantify the stress-free shape of the added particles; and (2) $\bfF_e := \bfF \bfF_{relax}^{-1}$ to quantify the elastic deformation, with $\bfsigma^* = \hat\bfsigma(\bfF_e)$, where $\hat\bfsigma(\cdot)$ is the stress response function of the added particles.
We will, in essence, define the reference configuration of the added particles through $\bfF_{relax}$ to ensure that we satisfy the kinematic compatibility of $\bfF$.

Substituting $\bfF = \bfF_e \bfF_{relax}$ in \eqref{eqn:def-grad-transport}, we find the equation governing the evolution of $\bfF_e(\bfx,t)$\footnote{
    We have assumed that there is no inelastic deformation in the material after it is attached to the body, which implies that $\bfF_{relax}$ is constant in time for each material particle, i.e., $\parderiv{\bfF_{relax}}{t} + (\bfv \cdot \nabla)\bfF_{relax}= 0$; \cite{naghibzadeh2021surface} discusses the more general case where $\bfF_{relax}$ can evolve in time.
}:
\begin{equation}
    \label{eqn:def-elastic-evolution}
    \parderiv{\bfF_e}{t}+(\bfv\cdot\nabla)\bfF_e=\left(\nabla\bfv \right) \bfF_e
\end{equation}
That is, $\bfF_e$ satisfies the same transport equation as $\bfF$, but with the important advantage that $\bfF_e$ has no kinematic compatibility requirements.
This makes the treatment of boundary conditions simple and direct, as opposed to grappling with issues of appropriate boundary conditions for \eqref{eqn:def-grad-transport} that preserve kinematic compatibility.
We therefore use $\bfF_e$ as the primary kinematic descriptor, and solve \eqref{eqn:def-elastic-evolution} instead of \eqref{eqn:def-grad-transport}.
The boundary condition at the growing boundary with $M>0$ is $\bfF_e = \hat\bfsigma^{-1}(\bfsigma^*)$\footnote{
    We assume that the stress response is invertible; however, from frame indifference, it can, at best, be invertible only up to a rotation. We can make any convenient choice for the unconstrained rotation.
}.
No boundary conditions on $\bfF_e$ are needed at the ablation boundaries ($M<0$).

\subsection{Constitutive Response}

To solve the balance laws in Section \ref{sec:balances}, we need to relate the stress in the growing body and the added material to the kinematic variables.
We use the framework of hyperelasticity, wherein the Cauchy stress tensor $\bfsigma$ only depends on the deformation gradient $\bfF$, and a scalar strain energy density function per unit mass $\hat W (\bfF)$ can be defined, such that:
\begin{equation}
    \bfsigma = \rho \parderiv{\hat W}{\bfF} \bfF^T
\end{equation}
Since the primary kinematic descriptor in this model is $\bfF_e$, we can redefine the elastic energy in terms of $\bfF_e$, and let $W(\bfF_e) = \hat W (\bfF)$, such that:
\begin{align}
     \label{eqn:hyperelasticMaterial}
    \bfsigma  = \rho \parderiv{W}{\bfF_e} \bfF_e^T
\end{align}

In all the problems considered in this paper, we use a compressible neo-Hookean material with the following strain energy function:
\begin{equation}
    \label{eqn: neo-hookean constitutive law}
    W(\bfF) = 
    \frac{\mu}{2} \left(\trace \left(\bfF^T \bfF \right) - \trace \left(\bfI\right) - 2 \log \det \left( \bfF \right) \right) + \frac{\lambda}{2}  \left(\det \left(\bfF \right)-1\right)^2
\end{equation}
where $\lambda$ and $\mu$ are the Lame parameters.
    We use the neo-Hookean energy because it provides a simple yet nonlinear energy function.
We have that $W(\bfI) = 0$ and $\left. \parderiv{W}{\bfF} \right|_{\bfF=\bfI} = \bfzero$.

\subsection{Model Summary}
\label{sec:model-summary}

In summary, our formulation models surface growth in hyperelastic materials by solving the coupled balances of mass and linear momentum, and the transport of elastic deformation, together with the appropriate boundary conditions and constitutive response \eqref{eqn:hyperelasticMaterial}. 
We use $\Omega_s$ to denote the growing solid body; $\Gamma_g$ to denote the growing boundary of the solid with $M>0$; $\Gamma_r$ to denote the growing boundary of the solid with $M<0$;  $\Gamma_0$ to denote non-growing part of the boundary ($M=0$), where $\partial \Omega_s = \Gamma_g \cup \Gamma_r \cup \Gamma_0$.
Further, we use $\Gamma_u$ to denote the boundaries where the velocity $\bfu_0$ is specified; and $\Gamma_t$ to denote the boundaries where the traction $\bft_b$ is specified, with $\partial \Omega_s = \Gamma_u \cup \Gamma_t$.

The balance of mass and the corresponding boundary condition read:
\begin{align}
    \label{eqn: continuuty-general}
    & \dot{\rho} = - \rho \divergence (\bfv) 
    \text{ on } \Omega_s,
    \qquad \text{ with }
    \\
    \label{eqn: continuuty-general-bc}
    & \bfV_b \cdot \hat{\bfn} = \bfv \cdot \hat{\bfn}  + \frac{M}{\rho} =
    \bfv \cdot \hat \bfn + \bfv_g \cdot \hat \bfn
    \text{ on } \partial \Omega_s
\end{align}

The balance of momentum and the corresponding boundary condition read:
\begin{align}
    \label{eqn: momentum-general}
    & \dot{\bfv} = \frac{1}{\rho} \divergence \bfsigma + \bfb
    \text{ on } \Omega_s 
    \\
    \label{eqn: momentum-general-bc}
    & \bfv = \bfv_0 \text{ on } \Gamma_u, \quad
    \bfsigma \hat{\bfn} 
    =
    M (\bfv_a - \bfv) + \bft_b
    \text{ on } \Gamma_t    
\end{align}

The transport of elastic deformation and the corresponding inlet boundary conditions read:
\begin{align}
    & \label{eqn: Fetransport-general}
    \dot{\bfF_e} = (\nabla \bfv) \bfF_e
    \text{ on } \Omega_s 
    \\
    \label{eqn: Fetransport-general-bc}
    & \bfF_e = \hat\bfsigma^{-1} (\bfsigma^*)
    \text{ on } \Gamma_g
\end{align}

We have used superposed dots to denote the material time derivative, i.e., $\dot{a} := \parderiv{a}{t} + \bfv \cdot \nabla a$.

\begin{remark}[Quasistatic simplification]
\label{sec:Quasi-Formulation}
In several surface growth processes like additive manufacturing and solidification in glaciers, material addition is very slow and inertial forces are negligible compared to the stress. 
Therefore the inertial terms in the momentum equation, i.e. $ \derivM{\bfv}{t}$ in \eqref{eqn: momentum-general} and $M\left( \bfv_a  - \bfv \right)$ in \eqref{eqn: momentum-general-bc}, can be neglected.
\end{remark}

\section{Numerical Method}
\label{sec:NumericalAlorithm}

In this section, we present the formulation of the numerical method.
Specifically, given the equilibrium density, elastic deformation, and the spatial domain occupied by the body at some discrete timestep, we aim to update these quantities, under the action of growth and mechanics, to their values at the next timestep.

Our method is implemented in Firedrake \cite{FiredrakeUserManual}, an open-source FE library, and the link to the code is provided at the end of the paper.
We use a continuous Galerkin polynomial interpolation for the displacement, elastic deformation, density, and phase function.
Although only 2-d examples are considered in this paper, the approach is readily extendable to 3-d.

This section is organized as follows:
\begin{itemize}

    \item 
    In Section \ref{sec:prelim-defns}, we provide definitions of the key quantities.
    
    \item 
    In Section \ref{sec:phaseFun}, we introduce a phase indicator function $\phi$ that implicitly tracks the shape of the growing solid body, and discuss its evolution and an associated regularization scheme.

    \item 
    In Section \ref{sec:incremental-quasistatics}, we discuss the incremental solution of the quasistatic hyperelastic problem that represents the balance of momentum.

    \item 
    In Section \ref{sec:balance-law-time}, we discuss the evolution equations for mass \eqref{eqn: continuuty-general}-\eqref{eqn: continuuty-general-bc}  and elastic deformation \eqref{eqn: Fetransport-general}-\eqref{eqn: Fetransport-general-bc}.
    The overall approach is to split the equations into source and transport contributions and solve them separately \cite{okazawa2007eulerian}.

    \item 
    In Section \ref{sec:growth-displacement-BC}, we discuss the time-stepping scheme that accounts for both growth and mechanics. 

    \item In Section \ref{sec:summary-numerical}, we summarize the overall method that brings together the components above.

\end{itemize}

\subsection{Preliminary Definitions}
\label{sec:prelim-defns}

The numerical method aims to evolve from one discrete time-step $t=t^n$ to the next time-step $t=t^{n+1} := t^n + \Delta t$.
Specifically, we are given the shape of the body, the spatial density distribution $\rho$, and the spatial elastic deformation distribution at equilibrium $\bfF_e$, and we describe how to find these quantities at the next time-step.
To perform this time-stepping, we decompose the evolution such that we solve \emph{separately} for growth and mechanical equilibrium.
That is, we consider the evolution of these quantities driven exclusively by growth and driven exclusively by mechanics (Fig. \ref{fig:mycaption-a}).
When we consider evolution driven by growth, we focus exclusively on the accretion and ablation of material particles without considering that the system is not in mechanical equilibrium; and when we focus on mechanics, we keep the set of material particles fixed and solve for mechanical equilibrium using the balance of momentum.

This leads to the following definitions:
\begin{itemize}
    \item 
    We denote by $\Omega^r$ the reference configuration; while we never construct the reference numerically, it is useful for various closed-form calculations.
    We denote the location of material particles in the reference by $\bfX$.

    \item
    We use $\bfx$ to denote spatial locations.

    \item 
    We denote by $\Omega^n$ and $\Omega^{n+1}$ the spatial domains occupied by the solid body at equilibrium at times $t^n$ and $t^{n+1}$.
    The spatial locations of the material particle located at $\bfX$ are $\bfx^{n}$ and $\bfx^{n+1}$ respectively.
    
    Note that $\Omega^{n+1}$ and $\Omega^n$ are the regions occupied by the equilibrium configurations of the body at different times and they include different sets of material particles due to growth, so they are not simply the mapping of the body under the deformation.
    
    \item 
    We denote by $\Omega^n_g$ the non-equilibrium intermediate configuration that arises from $\Omega^n$ after growth has occurred but before mechanical equilibrium occurs.
    The spatial location of the material particle located at $\bfX$ is $\bfx^{n}_{g}$.
    
    \item
    Traction boundary conditions are applied over the subset of the boundary denoted by $\Gamma_t^r$ and $\Gamma_{t_g}^n$ for the domains $\Omega^r$ and $\Omega_g^n$ respectively.
    
    \item
    We are focused here on a single timestep $\Delta t$ and assume that the velocity $\bfv$ is approximately constant over this timestep.
    Hence, the incremental displacement $\bfu:=\bfx^{n+1} - \bfx^{n}_{g}$ is related to the velocity by $\bfu \approx \bfv \Delta t$.

    \item
    Since we use a discretization that is fixed in space as the solid body deforms, we use a computational domain --- that we denote $R$ --- that contains the entire body as well as some part of the surrounding medium, such that the growing regions are always contained with $R$.

\end{itemize}
Given $\Omega^n$, $\rho^n(\bfx) := \rho(\bfx, t=t^n)$, and $\bfF_{e}^n(\bfx) := \bfF_{e}(\bfx, t=t^n)$, we develop a numerical algorithm that provides $\Omega^{n+1}$, $\rho^{n+1}(\bfx) := \rho(\bfx, t=t^{n+1})$, and $\bfF_e^{n+1}(\bfx) := \bfF_e(\bfx, t=t^{n+1})$.

We define extensions of various quantities from the growing solid body to the entire domain $R$, for reasons of numerical convenience.
The phase-field $\phi$ enables us to distinguish between the physically-relevant portion that is within the solid body, and the physically-irrelevant extensions.
For instance, we use a highly compliant linear elastic solid outside the physical body to provide smooth extensions of the displacement and mechanical velocity, and we use $\phi$ to appropriately construct the position-dependent elastic energy density.

\subsection{Phase Function} 
\label{sec:phaseFun}

We introduce a phase field $\phi$ that implicitly tracks the shape of the body, which is required to be able to apply the appropriate boundary conditions.
It is defined on $R$ such that $\phi(\bfx) \approx 1$ if $\bfx$ is within the solid body and $\phi(\bfx) \approx -1$ if $\bfx$ is outside it, i.e., at time $t=t^n$:
\begin{equation}
\label{eqn:phi-def}
    \phi(\bfx) \approx 
    \begin{cases}
    1 & \text{if }\bfx \in \Omega^n 
    \\
    -1 & \text{if }\bfx \in R \setminus \Omega^n
    \end{cases}
\end{equation}
For numerical convenience, $\phi$ is defined as a smooth function transitioning from $-1$ to $1$, where the level set $\phi=0$ represents the boundary of the region occupied by the growing body. 
Further, the outward unit normal of the solid body is $\hat \bfn = -\frac{\nabla \phi}{| \nabla \phi |}$.

For the evolution of $\phi$, we use a phase-field formulation that allows for arbitrary kinetics, following \cite{alber2005solutions,agrawal2015dynamic,agrawal2015dynamic-2,chua2022phase,guin2023phase,chua2024interplay}:
\begin{equation}
    \label{eqn:phaseFieldFull}
    \parderiv{\phi}{t} + \bfv_\phi \cdot \nabla \phi = 0 
\end{equation}
which represents the transport of level sets of $\phi$ with velocity $\bfv_\phi$, as discussed in \cite{agrawal2015dynamic}.

We highlight some important features of \eqref{eqn:phaseFieldFull}.
First, the only level set of physical significance is that corresponding to $\phi=0$, which represents the boundary of the solid body.
Second, $\nabla\phi$ is parallel to $\hat\bfn$, and hence only the normal component of $\bfv_\phi$ is of significance.
Hence, we define $\bfv_\phi$ such that it satisfies $\bfv_\phi=\bfV_b$ on the level set $\phi=0$, but are free to define any smooth extension of $\bfv_\phi$ over $R$.
We can then write \eqref{eqn: continuuty-general-bc} as $\bfv_\phi = \bfv + \bfv_g$ on the level set $\phi=0$.
For $\bfv_g$, we use a smooth extension of the growth velocity to $R$, with the specific choice of extension being chosen case-by-case for each problem.
For $\bfv$, which is the mechanical velocity, it is already defined over the solid body, and is extended over $R$ as described in Section \ref{sec:elastic-extension}.

We can write \eqref{eqn:phaseFieldFull} as:
\begin{equation}
    \label{eqn:phaseFieldFullSeparated}
    \parderiv{\phi}{t} + ( \bfv_g +  \bfv) \cdot \nabla \phi = 0 
\end{equation}
As mentioned in \ref{sec:prelim-defns}, at each timestep we solve \emph{separately} for mechanical equilibrium and growth.
Solving for mechanical equilibrium provides $\bfv$, and $\bfv_g$ is specified through additional growth physics.
We use a split-step method to evolve \eqref{eqn:phaseFieldFullSeparated} in time, i.e., we use each velocity individually to propagate a partial time step.
The specifics of solving this transport equation follow Section \ref{sec:Advection}.

\subsubsection{Regularization of the Phase Function}
\label{sec:reg}

The thickness of the regularized boundary of the solid body, i.e. the transition region where $\phi$ transitions from $-1$ to $+1$, is important for the numerical results.
If the thickness is smaller than the mesh size, it is not properly captured, while if it is very thick, the physical problem is not well represented.
Further, even if we start with $\phi$ having a moderate thickness, the evolution equation \eqref{eqn:phaseFieldFull} will often flatten or sharpen $\phi$ over time \cite{osher2005level}.
Finally, it has been observed that $\phi$ can potentially lose smoothness and regularity over time with the appearance of growing artificial oscillations.

We use the following approach to ensure the smoothness of $\phi$ over time during numerical calculations \cite{chen1992generation}. 
Let $\phi$ denote the numerical solution of \eqref{eqn:phaseFieldFullSeparated} at a given time step.
We apply the following variational filter to suppress potential oscillations and also prevent the interface from overly widening or sharpening:
\begin{equation}
    \label{eqn:reg_minimization}
    \mathrm{Reg}[\bar\phi] = 
    \int_R \left(\sigma_\phi \left( \bar\phi - \phi \right)^2 +
    \frac{\epsilon}{2} \left | \nabla \bar\phi \right |^2 + \frac{1}{2 \epsilon} \left( \bar\phi^2-1 \right)^2\right)
\end{equation}
Minimizing this functional provides a map $\phi \mapsto \bar\phi$ with the following features: (1) the first term controls the distance between $\bar\phi$ close to the physically-governed $\phi$; (2) the second term controls the interface smoothness to avoid spurious numerical oscillations; and (3) the third term controls the interface thickness.
$\sigma_\phi$ and $\epsilon$ are regularization parameters that weight the terms \cite{bronsard1991motion}.

The numerical values of the regularization parameters are problem-dependent, and are listed in the description of the problems that we solve in later sections.

\subsection{Momentum Balance through Quasistatic Incremental Hyperelasticity}
\label{sec:incremental-quasistatics}

We describe here how to find the displacement 
$\bfu(\bfx) = \bfx^{n+1} - \bfx^{n}_{g}$
to satisfy mechanical equilibrium, keeping fixed the set of material particles, i.e., without growth.

\subsubsection{Kinematics}

We begin by noting the useful identity:
\begin{equation}
\begin{split}
    & \bfx^{n+1}(\bfX) = \bfx^{n}_{g}(\bfX) + \bfu \left(\bfx^{n}_{g} \left(\bfX \right)\right) 
    \\ 
    \implies
	& 
    \parderiv{\bfx^{n+1}}{\bfX} = 
    \parderiv{\bfx^{n}_{g}}{\bfX} + \parderiv{\bfu}{\bfX} = 
    \parderiv{\bfx^{n}_{g}}{\bfX} + \parderiv{\bfu}{\bfx^{n}_{g}} \parderiv{\bfx^{n}_{g}}{\bfX}
    \\
	\implies
    &
    \label{eqn:evolutionF-Quasi}
    \bfF_e^{n+1}\bfF_{relax}^{n+1} 
    = 
    \bfF_{e_g}^n \bfF^n_{relax_g} + \left(\nabla_{\bfx^{n}_{g}} \bfu \right) \bfF_{e_g}^n \bfF^n_{relax_g}
    = 
    \left(\bfI + \nabla_{\bfx^{n}_{g}} \bfu \right) \bfF_{e_g}^n \bfF^n_{relax_g}
\end{split}
\end{equation} 
where $\nabla_{\bfx^{n}_{g}}$ is gradient with respect to ${\bfx^{n}_{g}}$. 
Using that $\bfF_{relax}$ does not evolve for a given material particle, it can be cancelled out to give:
\begin{equation}
    \bfF_e^{n+1}
    = 
    \left(\bfI + \nabla_{\bfx^{n}_{g}} \bfu \right) \bfF_{e_g}^n
\end{equation}

\subsubsection{Incremental Variational Principle}

Given the definition of $\bfsigma$ for a hyperelastic material \eqref{eqn:hyperelasticMaterial}, the solution of the balance of momentum \eqref{eqn: momentum-general} and \eqref{eqn: momentum-general-bc} in a quasistatic setting is equivalent to the minimization of the energy functional:
\begin{equation}
    I\left[\bfx^{n+1}\right] = 
      \int_{\Omega^r } \hat\rho_r(\bfX) \left(W\left(\bfF_e \left(\bfx^{n+1}\right)\right) - \bfb \cdot \bfx^{n+1} \right) \dm \bfX - 
     \int_{\Gamma_t^r } \bft_{b_r}\cdot \bfx^{n+1} \dm \bfX
\end{equation}
where $\bfb$, $\bft_{b_r}$, and $\hat\rho_r$ are the body force per unit mass, the traction at the boundary, and the density in the reference configuration $\Omega^r$.

We transform this to the configuration $\Omega_g^n$:
\begin{align}
    &
    I\left[\bfx^{n+1}\right] = 
    \int_{\Omega_g^n } \tilde\rho_r \left(\bfx_g^{n}\right)
    \left(W \left(\bfF_e \left(\bfx^{n+1} \right) \right) - \bfb \cdot \bfx^{n+1} \right) \det\left ( \parderiv{\bfx^{n}_{g}}{\bfX} \right)^{-1} \dm\bfx^{n}_{g} - 
    \int_{\Gamma_{t_g}^n} \bft_b \cdot \bfx^{n+1} \dm \bfx^{n}_{g}
    \label{eqn:integralxnX}
\end{align}
where $\bft_b$ is the traction in the configuration $\Omega_g^n$ and equal to $\mu_J^{-1} \bft_{b_r}$ with the area stretch $\mu_J$ between the configurations $\Omega^r$ and $\Omega_g^n$ \cite{gurtin1982introduction}, and $\rho_r = \hat{\rho}_r(\bfX) = \tilde\rho_r \left(\bfx^{n}_{g}\right)$ is the referential density.

We rewrite the energy functional \eqref{eqn:integralxnX} with respect to $\bfu$ instead of $\bfx^{n+1}$:
\begin{align}
    \label{eqn:I(u)}
    I[\bfu] = 
    \int_{\Omega_g^n } 
    \rho^n_g(\bfx^{n}_{g})
    \left(W \left( \left(\bfI + \nabla_{\bfx^{n}_{g}} \bfu \right) \bfF_{e_g}^n\left(\bfx^{n}_{g}\right) \right) - \bfb \cdot \left(\bfx^{n}_{g}+\bfu\right) \right)  \dm\bfx^{n}_{g}  - 
    \int_{\Gamma_{t_g}^n} \bft_b \cdot \left(\bfx^{n}_{g}+\bfu\right) \dm \bfx^{n}_{g}
\end{align}
where we have defined $\rho^n_g(\bfx^{n}_{g}) := \frac{\tilde\rho_r\left(\bfx^{n}_{g}\right)}{\det \bfF_{g}^n \left(\bfx^{n}_{g}\right)}$.

When we use small increments in time, we have that $\bfu$ is small, and approximate $\nabla_{\bfx^{n}_{g}}$ by $\nabla$, where $\nabla$ is the spatial gradient. 
Therefore, to find $\bfu$, we can minimize the incremental energy functional:
\begin{align}
    \label{eqn:I(u)Spatial}
    I[\bfu] = 
    \int_{\Omega_g^n } 
    \rho_g^n (\bfx)
    \left(W \left( \left(\bfI + \nabla \bfu \right) \bfF_{e_g}^n\left(\bfx\right) \right) - \bfb \cdot \left(\bfx+\bfu\right) \right)  \dm\bfx  - 
    \int_{\Gamma_{t_g}^n} \bft_b \cdot \left(\bfx+\bfu\right) \dm \bfx
\end{align}

\subsubsection{Extension to the Entire Computational Domain}
\label{sec:elastic-extension}

A simple, yet effective, approach to extend $\bfv$ beyond the solid body to the entire computational domain $R$ is to consider that a very soft isotropic linear elastic solid occupies the rest of the computational domain $R \backslash \Omega_s$, with elastic energy density function per unit mass given by:
\begin{equation}
\label{eqn:lineaElasticity}
    W_\text{linear} \left( \nabla\bfu \right) = \frac{\mu_c}{4} 
           \left|\nabla \bfu + \left( \nabla \bfu \right)^T\right|^2 + 
           \frac{\lambda_c}{2} \trace \left(\nabla \bfu \right) 
\end{equation}
where $\mu_c$ and $\lambda_c$ are the Lame parameters with numerical values smaller by a factor $\approx 10^{-3}$ than the physical solid.
From our numerical examples, we find that these values are sufficiently small that any numerical artifacts are negligible.

We solve \eqref{eqn:I(u)Spatial} over $R$ instead of $\Omega_g^n$, and use an elastic energy density of the form:
\begin{equation}
    \label{eqn:WFull}
    W \left( \left( \bfI + \nabla \bfu \right) \bfF_{e_{g}}^n \right) 
    = 
        H_l(\phi(\bfx))  W_s \left( \left(\bfI + \nabla \bfu \right) \bfF_{e_{g}}^n \right) 
        + 
        (1-H_l(\phi(\bfx))) W_\text{linear} \left(\nabla \bfu \right) 
\end{equation}
where $W_s$ is the elastic energy of the solid, and $H_l(\phi)$ is a function that transitions rapidly and smoothly from $0$ to $1$ when the argument changes sign, such that $H_l(\phi=1) = 1$ and $H_l(\phi=-1) = 0$.
Throughout this paper, we use $H_l(\phi) = \frac{1 + \tanh (l \phi) }{2}$ where $l$ is a modeling constant prescribing the width of the transition zone.

\subsubsection{Boundary Conditions}
\label{sec:BC}

For boundaries with prescribed displacement that is constant in time, we use a computational domain such that those boundaries of the solid body coincide with the boundaries of the computational domain, and we treat them as in standard FE approaches.
While we do not consider problems in which the prescribed boundary evolves, the diffuse Nitsche method provides a promising approach to impose the displacement on the boundary \cite{nguyen2018diffuse}.

For boundaries with prescribed traction $\bft_b$, we use a computational domain such that these boundaries lie within the domain. 
We then smear out the traction vector over the width of the interface and replace it with a body force equal to $\bft_b {|\nabla H_l(\phi)|}$ using the following approximation \cite{nguyen2018diffuse}:
\begin{equation}
\label{eqn:tractionBC}
    \int_{\Gamma_{t_g}^n} \bft_b \cdot \left(\bfu+\bfx^{n}_{g} \right) \dm \bfx^{n}_{g} = 
    \int_R \bft_b \cdot \left( \bfu+\bfx^{n}_{g} \right) \delta_{\Gamma_{t_g}^n} \dm V
    \approx
    \int_R \bft_b \cdot \left(\bfu+\bfx^{n}_{g} \right) {|\nabla H_l(\phi)|} \dm V
\end{equation}
where $\delta_{\Gamma_{t_g}^n}$ is a Dirac mass supported on the boundary $\Gamma_{t_g}^n$.

\subsection{Evolution of Elastic Deformation and Density}
\label{sec:balance-law-time}

The evolution equations for density \eqref{eqn: continuuty-general} and elastic deformation \eqref{eqn: Fetransport-general-bc} have the common structure:
\begin{equation}
    \parderiv{f}{t} + \bfv \cdot \nabla f = \cal{F}
\end{equation}
where $f$ can be $\rho$ or $\bfF_e$. 
Using a finite-difference approximation to the time derivative, we can split the equations to separate the source and transport contributions, and then solve these separately \cite{okazawa2007eulerian}.

\subsubsection{Transport Contribution}
\label{sec:Advection}

We first turn to the transport contribution that is discretized in time using an explicit difference to give:
\begin{equation}
    \label{eqn: advection-general}
    \frac{f^{n+1} - f^{n}_{g}}{\Delta t} + \bfv \cdot \nabla  f^{n}_{g} = 0
\end{equation}

Let $\cal T$ denote a discretization of $R$, and $\bar\bfx_i$ denote the spatial location of the grid point indexed by $i$ in $\cal T$. 
For a fixed discretization, $\bar\bfx_i$ remains constant in time. 
We denote the nodal values $f^{n}_{g}(\bar\bfx_i)$ by $f_{g_i}^n$ and $f^{n+1}(\bar\bfx_i)$ by $f_i^{n+1}$.
Here, we describe how to evaluate $\left\{f\right\}^{n+1}$ given $\left\{f\right\}^n_g$ and $\bfv (\bar\bfx_i)$.

Recall that $\bfx^{n}_{g}$ and $\bfx^{n+1}$ are the spatial locations of the material particle $\bfX$ in the configurations $\Omega^n_g$ and $\Omega^{n+1}$.
Hence, we have:
\begin{equation}
    f^{n+1}\left(\bfx^{n+1}\right) = f^{n}_{g}\left(\bfx^{n}_{g}\right) 
\end{equation}
Using $\bfx^{n}_{g} \approx \bfx^{n+1} - \bfv\left(\bfx^{n+1}\right)\Delta t$ in the argument of the left side in the equation above gives:
\begin{equation}
    f^{n+1}\left(\bfx^{n+1}\right) = f^{n}_{g}\left(\bfx^{n+1} - \bfv\left(\bfx^{n+1}\right)\Delta t\right) 
\end{equation}
Setting $\bfx^{n+1} = \bar\bfx_i$ gives:
\begin{equation}
\label{eqn:transport-explicit}
    f_i^{n+1} = f^{n}_{g}\left(\bar\bfx_i - \bfv\left(\bar\bfx_i\right)\Delta t\right) 
\end{equation}
Since $\bar\bfx_i -  \bfv(\bar\bfx_i) \Delta t$ does not generally lie on any mesh point, we use interpolation to evaluate $f^{n}_{g} \left(\bar\bfx_i - \bfv(\bar\bfx_i) \Delta t \right)$ inside an element containing $\bar\bfx_i - \bfv(\bar\bfx_i) \Delta t$. 
The case when $\bar\bfx_i - \bfv(\bar\bfx_i) \Delta t$ is outside $R$ is discussed in the examples below.

We note that an explicit approach does not require the solution of a system of nonlinear algebraic equations and is computationally efficient while maintaining stability and accuracy, as we see in the examples.

\subsubsection{Source Contribution}
\label{sec:source}

\paragraph{Elastic Deformation.}

We use implicit Euler time discretization to write:
\begin{equation}
\label{eqn:Fe-increment}
    \frac{\bfF_e^{n+1} (\bfx) - \bfF_{e_g}^n(\bfx)}{\Delta t} = \left(\nabla \bfv\right)  \bfF_e^{n+1}(\bfx) 
    \implies
    (\bfI - \nabla \bfu)  \bfF_e^{n+1}(\bfx) = \bfF_{e_g}^n(\bfx)
\end{equation}
For small $\nabla \bfu$, we can use the approximation $\left(\bfI - \nabla \bfu \right)^{-1} \approx \left( \bfI + \nabla \bfu \right)$ to write:
\begin{equation}
    \label{eqn: eqn:Fe-increment-Final}
     \bfF_e^{n+1} (\bfx) = 
    (\bfI + \nabla \bfu) \bfF_{e_g}^{n} (\bfx).
\end{equation}

\paragraph{Density.}

We use implicit Euler time discretization to write:
\begin{equation}
\label{eqn:Continuuty-increment}
    \frac{ \rho^{n+1} (\bfx) - \rho_g^{n}\left(\bfx \right)}{\Delta t} = 
    - \left( \divergence \bfv \right)  \rho^{n+1} \left(\bfx\right) 
    \implies 
    \left( 1 + \divergence \bfu \right) \rho^{n+1}\left(\bfx \right)  = 
    \rho_g^{n}(\bfx)
\end{equation}
For small $\nabla \bfu$, we can use the approximation $\left(1+\divergence \bfu \right) \approx \det \left( \bfI + \nabla \bfu \right)$ to write:
\begin{align}
    \label{eqn:rho_proof}
     \rho^{n+1}(\bfx) =  
    \frac{\rho_g^{n}(\bfx)}{\det \left( \bfI + \nabla \bfu\right)}
\end{align}

We note that although we use an implicit approach, our linearization enables us to bypass solving a system of algebraic equations, making the implicit scheme computationally efficient.

\subsubsection{Extension to the Entire Computational Domain}

To obtain a smooth extension of $\bfF_{e}$ and $\rho$ over $R$, we simply use that we have considered the entirety of $R$ to be occupied by a solid body following Section \ref{sec:elastic-extension}.
Consequently, all quantities associated with the physical solid body are also extended smoothly to the soft surrounding solid, using an interpolation of the form in \eqref{eqn:WFull}.

\subsection{Time-stepping with Growth and Mechanics}
\label{sec:growth-displacement-BC}

\begin{figure}[htb!]
    \subfloat[Unconstrained Growth]{{\includegraphics[height=.17\textwidth]{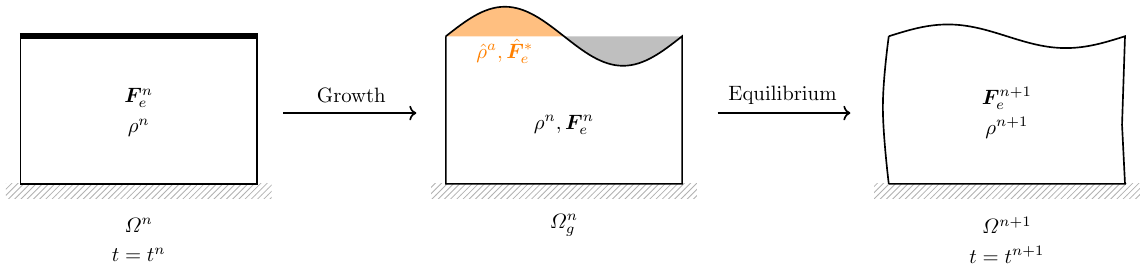} }\label{fig:mycaption-a}}
    \\
    \subfloat[Constrained Growth]{{\includegraphics[height=.17\textwidth]{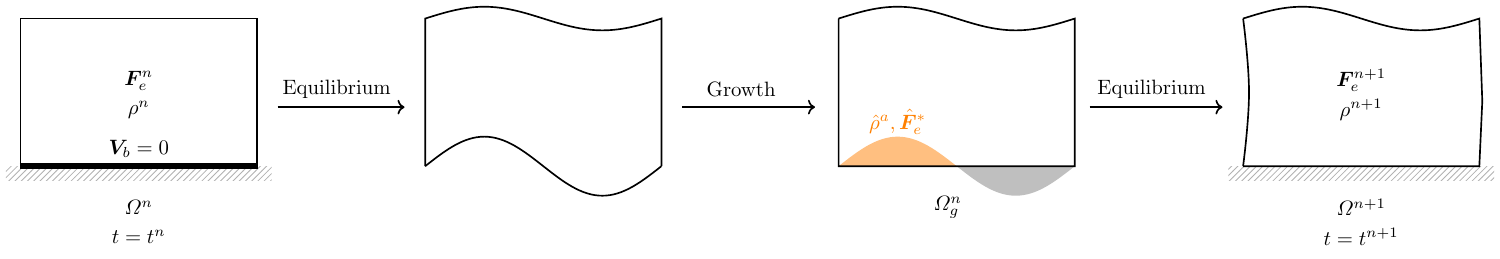} }\label{fig:mycaption-b}}
    \caption{A schematic showing the calculations within a single time-step. The thick line depicts the growing boundary with prescribed growth velocity $\bfv_g \cdot \hat \bfn$.  
    (a) Unconstrained Growth: growth occurs at a boundary where the displacement is not prescribed.
    (b) Constrained Growth: growth occurs at a boundary where the displacement is prescribed, thereby constraining the relation between the boundary velocity, growth velocity, and mechanical velocity; the figure shows the case that $\bfv_g$ and $\bfv$ are coupled through $\bfV_b = \bfzero$. 
    During growth, we do not consider mechanics, and during mechanical relaxation, we do not consider growth.
    In $\Omega_g^n$, the region that is not shaded consists of material from the previous time-step $t=t^n$, the region shaded orange is the accreted material with $\rho=\hat\rho_a, \bfF_e = \hat \bfF_e^*$, and the region shaded grey is the ablated material.}
    \label{fig:independent vs dependent}
\end{figure}

We discuss two different physical situations which must be distinguished in the numerical method.
With reference to Figure \ref{fig:independent vs dependent}, the top row shows growth occurring at a boundary at which the displacement is not prescribed, which we denote ``Unconstrained Growth'', and the bottom row shows growth occurring at a boundary at which the displacement is prescribed, which we denote ``Constrained Growth''.
In the notation of Section \ref{sec:model-summary}, Constrained Growth corresponds to the case where there is nontrivial $\Gamma_u \cap \Gamma_g$.
Specifically, if growth occurs at portions of the boundary where the displacement is prescribed, we must take care to ensure that both the specified growth velocity conditions and the displacement boundary conditions are satisfied.
We first discuss the case of Unconstrained Growth.

\subsubsection{Unconstrained Growth}
\label{sec:UG}

We recall the decomposition in \eqref{eqn: continuuty-general-bc} that $\bfV_b \cdot \hat \bfn = \bfv_g \cdot \hat \bfn + \bfv \cdot \hat \bfn$ on $\partial\Omega_s$, i.e., the motion of the boundary is the superposition of the growth velocity $\bfv_g$ and the material velocity $\bfv$.
The case of Unconstrained Growth is straightforward: when $\bfV_b$ is not prescribed, we can compute the shape of the body due to growth using the prescribed growth velocity $\bfv_g$, and then use mechanical equilibrium with this added material to find $\bfv$.

Within a single time-step, we perform the following calculations, referring to Figure \ref{fig:independent vs dependent}(a).

\begin{enumerate}

    \item
    We evolve $\phi$ using the growth velocity $\bfv_g$, which changes the shape of the body due to the ablation and accretion of material, but we do not consider mechanics.

    \item 
    We solve for mechanical equilibrium of the body with the ablated and accreted material to account for the deformation due to the added material.
    This provides the displacement of the body due to mechanical relaxation.

    \item
    We use the mechanical velocity $\bfv$ that is computed from mechanical relaxation in the previous step to transport $\phi$, $\rho$, $\bfF_e$.
    The evolution of $\phi$ in this step provides the change in shape of the body due to mechanical relaxation.
    
\end{enumerate}

\subsubsection{Constrained Growth}
\label{sec:CG}

In the case of Constrained Growth, both $\bfV_b$ and $\bfv_g$ are prescribed on some parts of the boundary.
Hence, $\bfv \cdot \hat \bfn = \bfV_b \cdot \hat \bfn - \bfv_g \cdot \hat \bfn$ is constrained.
To deal with this constraint, we perform the following calculations in a single time-step (Fig. \ref{fig:independent vs dependent}(b)).

\begin{enumerate}
    \item 
    We solve for mechanical equilibrium of the body without growth, using a displacement boundary condition that violates the prescribed displacement but ``makes room'' for the ablation and accretion of material.
    That is, we use the boundary condition corresponding to $\bfV_b - \bfv_g$ prescribed, rather than $\bfV_b$ prescribed, on the relevant part of the boundary when we solve for mechanical equilibrium.

    \item
    We use the mechanical velocity $\bfv$ that is computed from mechanical relaxation in the previous step to transport $\phi$, $\rho$, $\bfF_e$.
    The evolution of $\phi$ in this step provides the change in shape of the body due to mechanical relaxation.

    \item
    We ``fill in'' the ablation and accretion of material using $\bfv_g$, which restores the required displacement boundary condition that $\bfV_b$ is prescribed.

    \item 
    We solve for mechanical equilibrium of the body with the ablated and accreted material, with the displacement-prescribed boundary conditions,  to account for the deformation due to the added material.
    This provides the displacement of the body due to mechanical relaxation.

    \item
    We use the mechanical velocity $\bfv$ that is computed from mechanical relaxation in the previous step to transport $\phi$, $\rho$, $\bfF_e$.
    The evolution of $\phi$ in this step provides the change in shape of the body due to mechanical relaxation.

\end{enumerate}

\subsection{Summary}
\label{sec:summary-numerical}

We summarize the method in Algorithm \ref{alg:Growth-U}, showing the steps go from the state at $t=t^n$ to $t=t^{n+1}$.
We focus on Unconstrained Growth, but Constrained Growth is similar following the procedure in Section \ref{sec:CG}.

\begin{algorithm}[H]
\caption{Unconstrained surface growth of a deformable body, from $t=t^n$ to $t=t^{n+1}=t^n + \Delta t$. }
\label{alg:Growth-U}

\begin{algorithmic}[1]
\State \textbf{Initialize:} 
Given the initial geometry $\phi^n$ (representing $\Omega^n$), 
initial elastic deformation $\bfF_{e}^n(\bfx)$, 
initial density $\rho^n(\bfx)$ of the body, 
boundary conditions,
growth velocity $\bfv_g(\bfx)$,
and time step size $\Delta t$ at timestep $n$.

    \State \textbf{Compute the non-equilibrium growing body $\Omega_g^n$:} 
    $\phi_g^{n} \left(\bfx \right) \gets \phi^n \left(\bfx -  \bfv_g \Delta t \right)$
    following \eqref{eqn:transport-explicit} and the discussion below \eqref{eqn:phaseFieldFullSeparated}.

     \State \textbf{Construct  of $\bfF_{e_g}^n$ and $\rho_g^n$:} Assigning density $\hat \rho^a$ and elastic deformation $\hat \bfF_e^*$ of the added material to the accreted region. 
    
    \State \textbf{Relax to mechanical equilibrium:} Minimizing the energy functional \eqref{eqn:I(u)Spatial} over $R$ to find $\bfu$.
    
    \State \textbf{Update the elastic deformation and density based on source contribution:}
    $\tilde \bfF_e^{n+1}(\bfx) \gets \left(\bfI + \nabla \bfu\right) \bfF_{e_g}^n \left(\bfx\right)$ following \eqref{eqn: eqn:Fe-increment-Final},
    and $\tilde \rho^{n+1}(\bfx) \gets \frac{\rho_g^{n} \left(\bfx \right)}{\det \left( \bfI + \nabla \bfu\right)}$ following \eqref{eqn:rho_proof}.
    
    \State \textbf{Update the elastic deformation and density based on transport contribution:} 
    $\bfF_{e}^{n+1} \left(\bfx \right) \leftarrow \tilde \bfF_e^{n+1}\left(\bfx - \bfu \right)$ 
    and 
    $\rho^{n+1} \left( \bfx \right) \tilde \rho^{n+1} \left(\bfx - \bfu \right)$
    following \eqref{eqn:transport-explicit}.
    
    \State \textbf{Advect the phase-field to generate $\Omega^{n+1}$:} 
    $\tilde \phi \left(\bfx \right) \gets \phi_g^n \left(\bfx-\bfu \right)$ following \eqref{eqn:transport-explicit} and the discussion below \eqref{eqn:phaseFieldFullSeparated}.
    
    \State \textbf{Filter $\tilde \phi$ for regularity:} Minimizing \eqref{eqn:reg_minimization} given $\tilde \phi \left(\bfx \right) $ to evaluate $\phi^{n+1}$.
    
    \State \textbf{Update the initial condition for the next timestep:} Assign fields $\phi^{n+1}$, $\rho^{n+1}$ and $\bfF_{e}^{n+1}$ to $\phi^{n}$, $\rho^{n}$ and $\bfF_{e}^{n}$.
    
\end{algorithmic}
\end{algorithm}

\section{Non-normal Growth}
\label{sec:examples}
\label{sec:Non-normalGrowth}

In the growth of hard biological tissues such as horns and nails, growth occurs at a fixed surface, and the generating cells that push the growing material outward might not be aligned with the normal vector of the surface \cite{skalak1997kinematics}. 
Hence the growth of a stress-free material in a non-normal direction will be observed.
Although complex horn shapes can be produced using methods based solely on kinematics \cite{moulton2014surface, garikipati2009kinematics}, the effect of stress and deformation in non-normal growth has not been studied because of the challenges in an evolving reference configuration that grows in a non-normal direction.

In this section, we aim to simulate non-normal growth of stress-free objects against a fixed boundary.
We first highlight that in the phase field approach, only the normal component of the growth velocity can affect the evolution of the phase function, following the discussion in Section \ref{sec:phaseFun} and also be noticed from \eqref{eqn:mass-bc-main} and \eqref{eqn:mom-bc-main}.
To achieve non-normal growth, we model growth with a normal velocity but further add that the accreted material is not stress-free.
In essence, we will find that the relaxation of the accreted material leads to non-normal growth.

We begin with an example that computes the relaxation to mechanical equilibrium without growth to verify the approach for incremental hyperelasticity and the transport of density and elastic deformation.

\begin{figure}[htb!]
    \centering
    \includegraphics[width=0.95\textwidth]{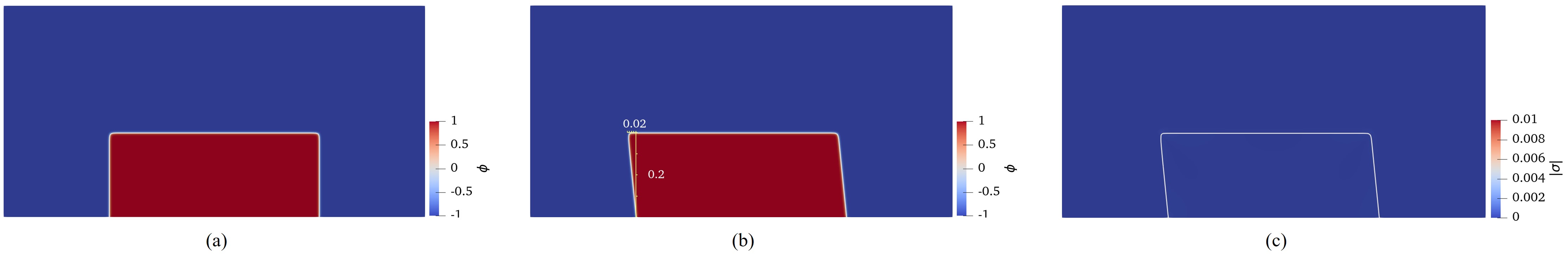}
    \caption{Deformation of a pre-stressed rectangular solid body: (a) the initial shape of the body, (b) the equilibrium shape of the body, and (c) the norm of the Cauchy stress tensor at equilibrium. The white line represents the $\phi=0$ levelset showing the boundary of the body.}
    \label{fig:ElasticityNonNormal}
\end{figure}

\begin{example}[Relaxation of a pre-stressed solid body without growth]
    We consider a rectangular solid body modeled by a neo-Hookean material \eqref{eqn: neo-hookean constitutive law}, which is minimized when the argument is the identity $\bfI$.
    We apply a pre-stress by setting the initial elastic deformation to have the value 
    $\bfF_e(\bfx,t=0) = \begin{bmatrix}
        1 & 0.1 \\
        0 & 1
        \end{bmatrix}$
    uniformly in the body.
    The lower face is clamped to a wall with displacement prescribed to be $\bfzero$ and the other faces are traction-free.
    
    Figure \ref{fig:ElasticityNonNormal} shows the initial rectangular shape and the relaxed shape after we solve for mechanical equilibrium as well as the stress-free equilibrium state.
    Furthermore, we highlight that the relatively large shear is beyond the linear regime and is captured well by the algorithm, with the final deformation corresponding to     $\begin{bmatrix}
        1 & 0.1 \\
        0 & 1
        \end{bmatrix}^{-1}$.  
\end{example}

We next study non-normal growth.
We use $\theta$ to denote the angle of the generating cells with the normal of the fixed boundary.
The accreted material is defined to have a pre-stress that corresponds to a shear strain of $\alpha = \tan{(-\theta)}$, i.e., $\hat \bfF_{e}^* = \begin{bmatrix}
    1 & \tan(-\theta) \\
    0 & 1
    \end{bmatrix}$, and the density $\hat \rho^a$ of the accreted material is set to $1$.
Accretion occurs through the addition of mass at a rate $\frac{M}{\rho} = v_0 > 0$ from the bottom face that is also clamped.
As discussed in Section \ref{sec:growth-displacement-BC}, this is a case of constrained growth with the growth velocity and boundary velocity both being prescribed on a common part of the boundary.

We note two computational details.
First, we do not need a smooth extension of the growth velocity since it is specified only at a part of the boundary of $R$; we simply use a growth velocity of $0$ throughout $R$.
Second, in the transport update (Section \ref{sec:Advection}), when we require information from outside the domain, this only occurs near the accretion boundary.
In that case, we use $\bfF_{e} = \begin{bmatrix}
    1 & \tan(-\theta) \\
    0 & 1
    \end{bmatrix}$, $\rho=1$, and $\phi=1$.

\begin{figure}[htb!]
    \subfloat[]{{\includegraphics[width=.3\textwidth]{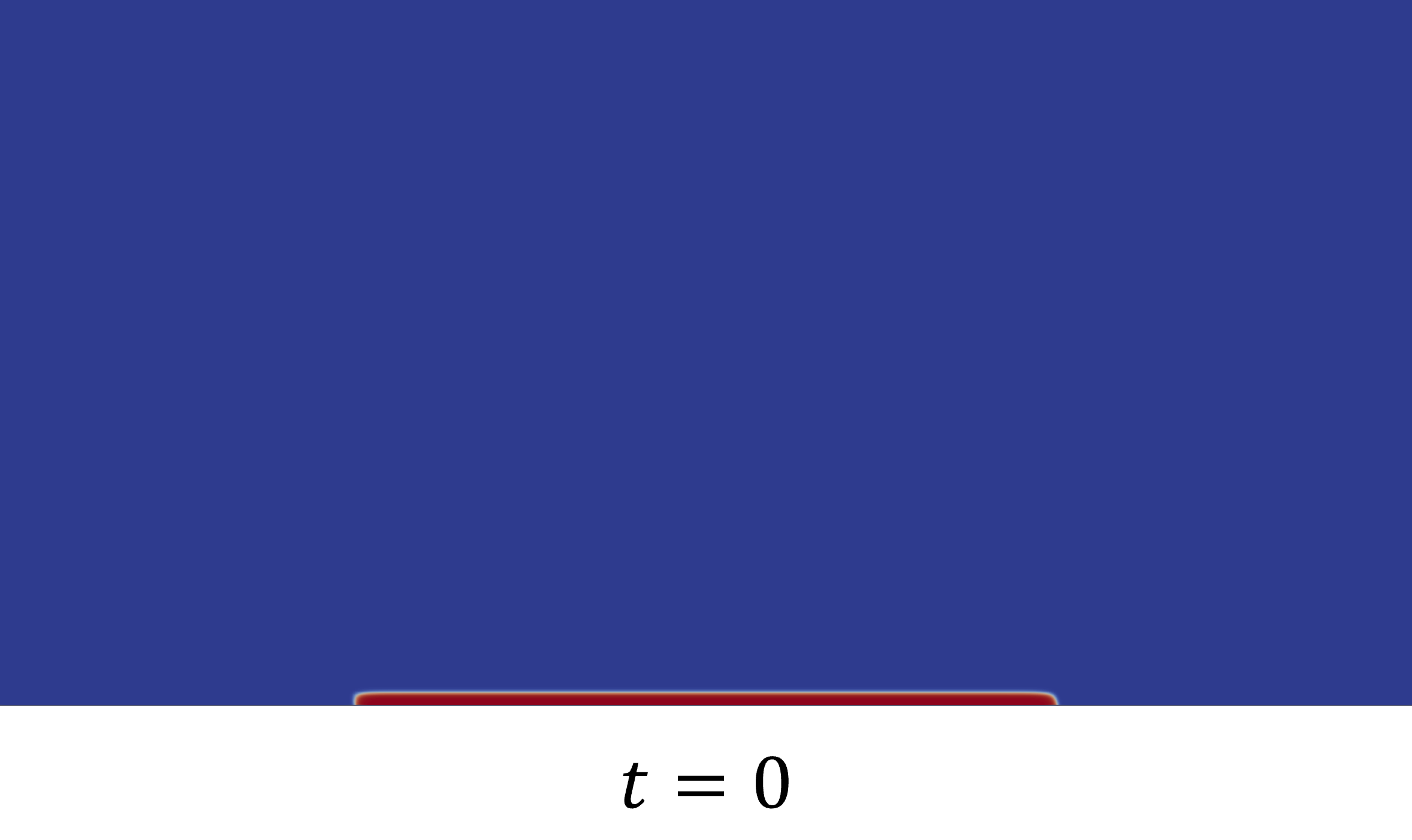} }\label{fig:constant-a1}
    {\includegraphics[width=.3\textwidth]{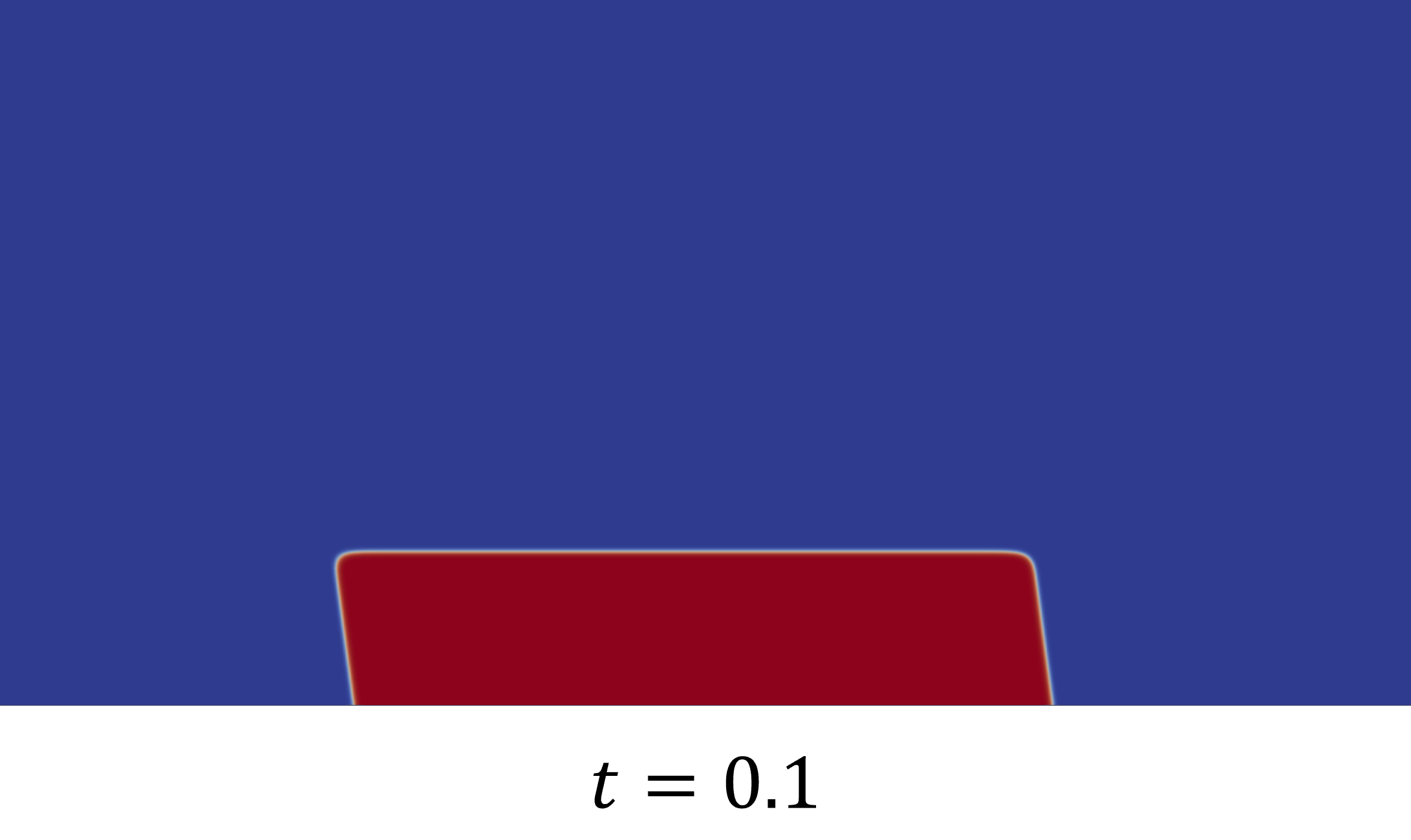} }\label{fig:constant-a2}
    {\includegraphics[width=.384\textwidth]{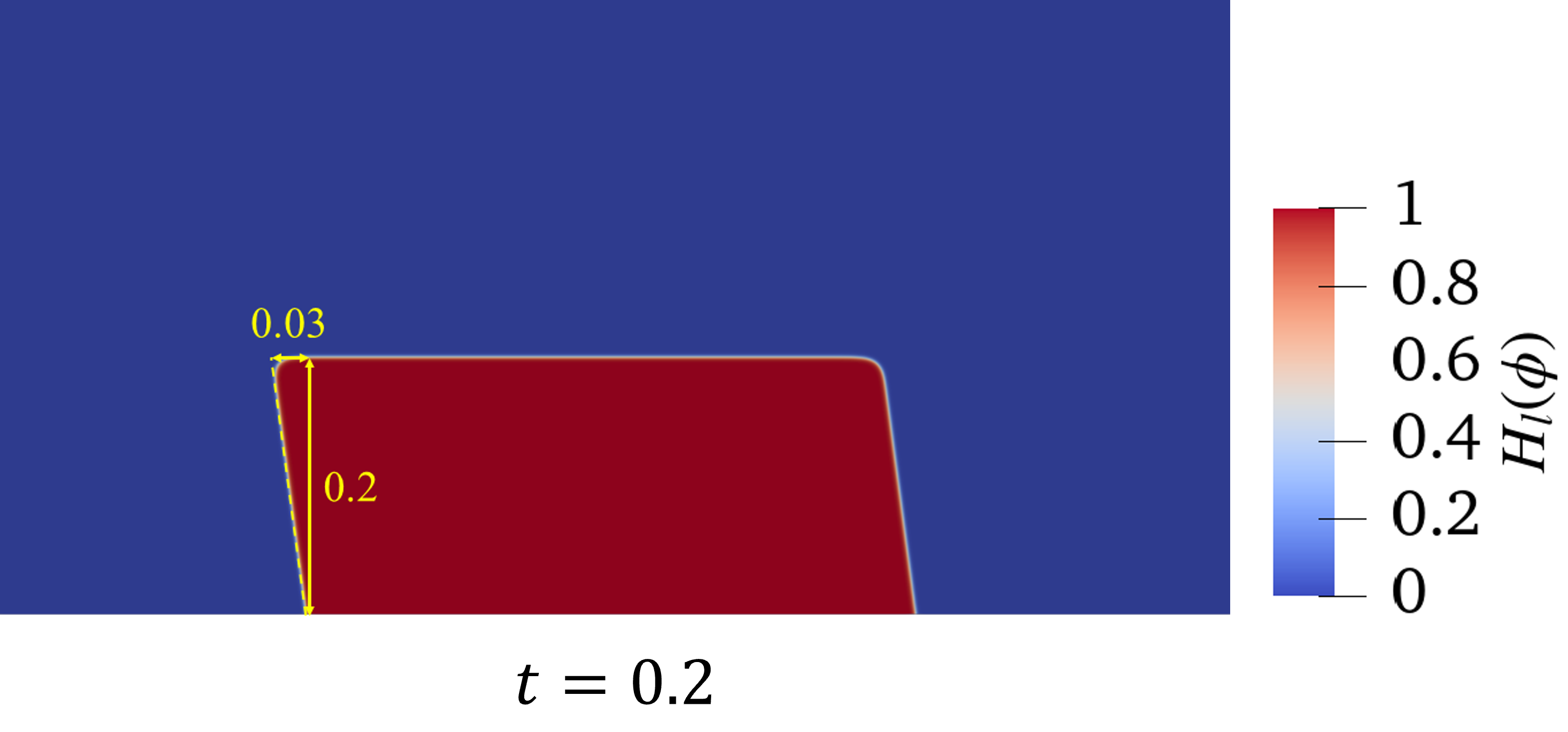} }\label{fig:constant-a3}}
    \\
    \subfloat[]{{\includegraphics[width=0.984\textwidth]{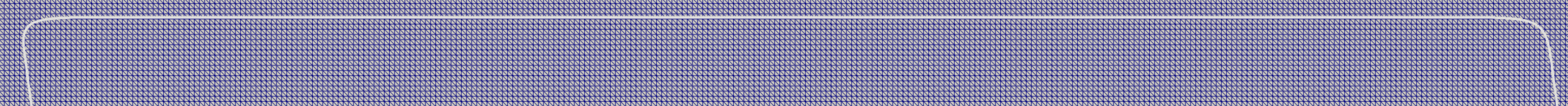} }\label{fig:constant-b}}
    \\
    \subfloat[]{{\includegraphics[width=.3\textwidth]{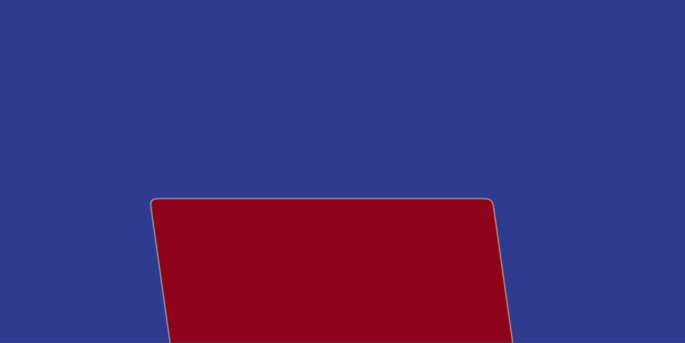} }\label{fig:constant-c}}
    \subfloat[]{{\includegraphics[width=.395\textwidth]{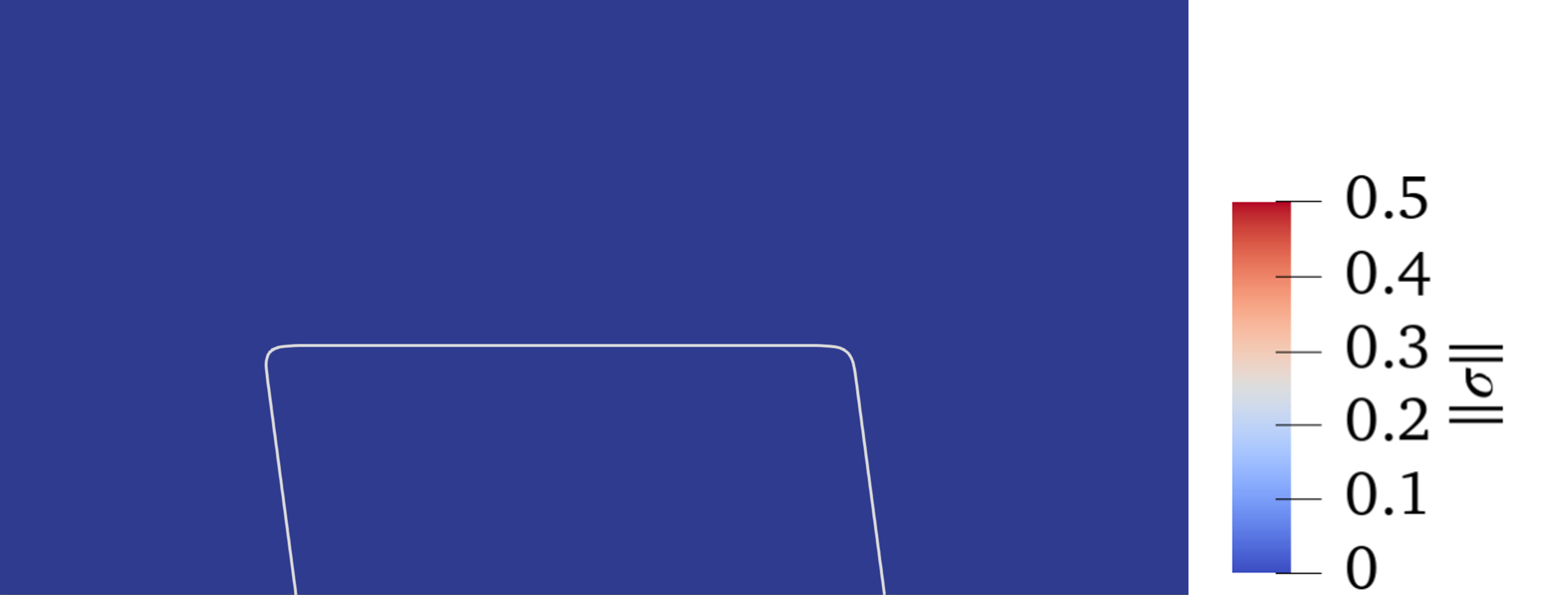} }\label{fig:constant-d}}
    \caption{Non-normal growth from a fixed surface. (a) The evolution of the phase function is shown at different times, (b) a zoomed-in view of the computational mesh shown at the end of the second time step, where the thick white line indicates the boundary of the growing body, (c) the result for a mesh that is finer by a factor of 2, at $t=0.2$, suggesting that the method is robust, and (d) the norm of the Cauchy stress tensor at $t=0.2$. The white curve shows the zero levelset of $\phi$ that indicates the shape of the growing body.}
    \label{fig:NNG}
\end{figure}

\begin{example}[Non normal growth with a fixed $\theta$]
    We consider $\alpha = \tan(\theta) = 0.15$ and $v_0 = 1$ and constant in time.  
    Figure \ref{fig:NNG} shows the non-normal growth of the growing solid body at different times.
    We observe from the levelset showing the shape of the body that there is non-normal growth and it occurs in the direction prescribed by our choice of $\alpha$.
    Further, the equilibrium stress distribution is essentially $\bfzero$ everywhere.
    
    We note that the computational domain has width $2$ and height $1$; the growth boundary has width $1$; the mesh size $\Delta x$ is $\approx 1.6 \times 10^{-3}$; $\epsilon = 3 \Delta x$; $\sigma_\phi = \frac{5}{\Delta x}$; and $\Delta t= 0.01$, in units that are non-dimensionalized by the growth velocity.
\end{example}

\begin{figure}[htb!]
    \centering
    \includegraphics[width=.98\textwidth]{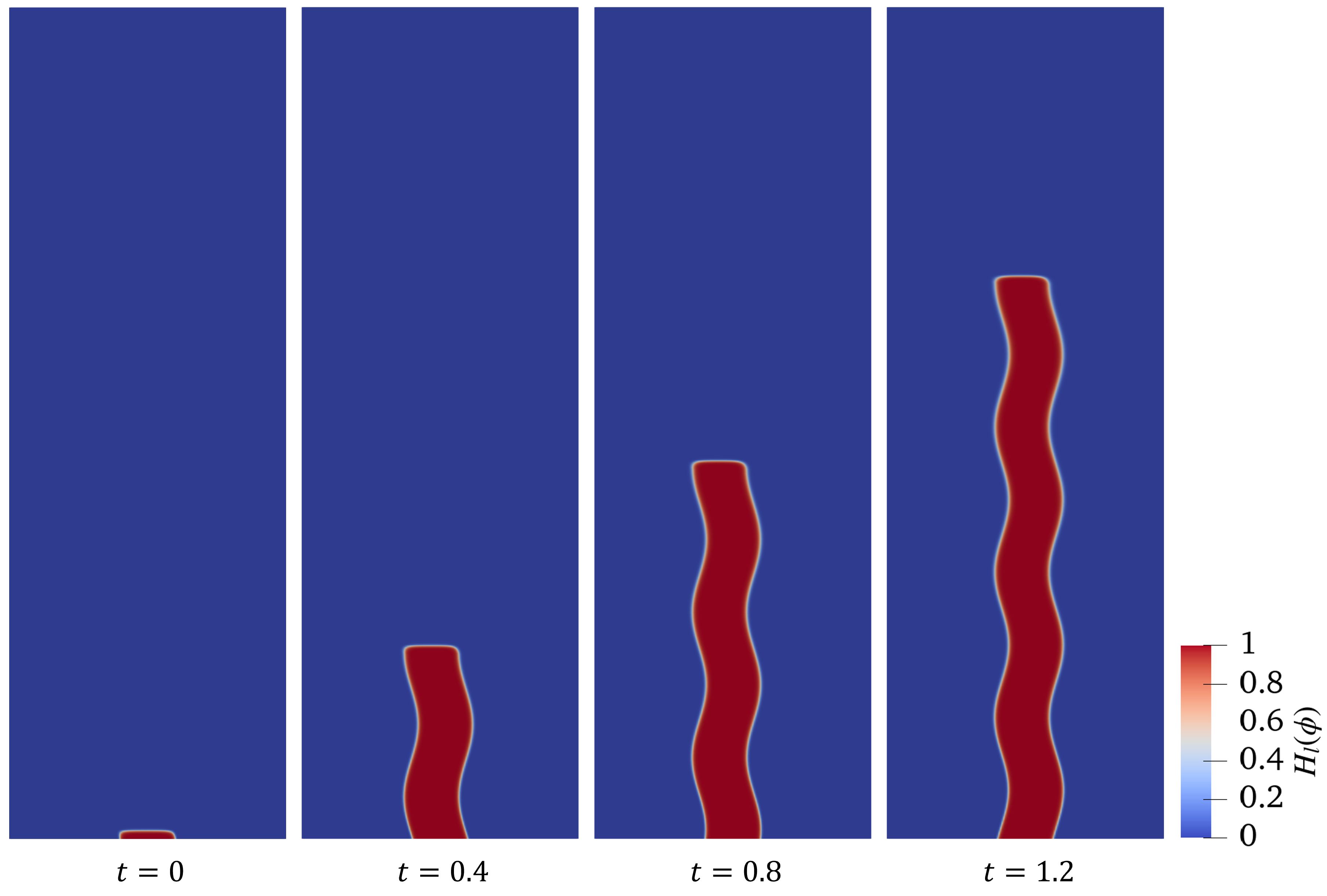}
    \caption{Non-normal growth from a fixed surface with a time-varying pre-stress corresponding to $\alpha = 0.3 \sin(20 t)$. The evolution of the phase function is shown at various times, and the white curve shows the zero levelset of $\phi$ that indicates the shape of the growing body.  }
    \label{fig:NNGOsci}
\end{figure}

\begin{example}[Non-normal growth with a varying $\theta$]
    We demonstrate the growth of complex shapes due to non-normal growth when $\theta$ varies in time: a growing body at different times using $\alpha = 0.3 \sin(20 t) $ and $v_0 = 0.5$ is shown in Figure \ref{fig:NNGOsci}.
    The computational domain has non-dimensional width of $0.3$ and height of $0.9$, and the growth surface has a width of $0.06$; $\Delta t= 0.01$; and the mesh has element size $\approx 1.6 \times 10^{-3}$ with the same regularization parameters as in the previous example.
\end{example}

    The examples in this section are formulated such that the shape of the stress-free body is known, enabling us to verify the accuracy of the numerical implementation.
    This also confirms that the soft isotropic linear material used in the exterior of the physical solid has a negligible effect.

\section{Application to Regelation: Stress-driven Melting and Refreezing}
\label{sec:MeltingUnderPressure}

Stress-induced melting occurs when a solid, at a constant temperature that is below the melting temperature, undergoes a phase transition into the liquid state due to the application of external load; once the load is removed, the liquid refreezes.
This phenomenon is called regelation and occurs when the solid state density is different from that of its liquid state.
A well-known illustration of this phenomenon is hexagonal ice --- the most abundant ice on the earth --- whose density is approximately $10\%$ lower than that of water, and which is often not too far below the melting temperature in environmental settings.
Understanding the dynamics of the stressed ice and water interface is important for various environmental mechanics issues \cite{dash1995premelting} such as glaciology \cite{weertman1957sliding} and frost heave analysis in soil, permafrost, and infrastructure \cite{konrad1980mechanistic}, as well as cryobiology \cite{taylor2019new}.

Unlike the model problems considered in Section \ref{sec:Non-normalGrowth}, here we do not arbitrarily specify the growth velocity and properties of the accreted material.
Instead, we obtain it using thermodynamics, by formulating a thermomechanical free energy and applying the Coleman-Noll procedure to construct an evolution equation that is consistent with the laws of thermomechanics \cite{gurtin1982introduction}\footnote{We assume that the temperature does not evolve in time, and hence do not need to solve the energy equation explicitly beyond its role in the Coleman-Noll analysis.}.
Specifically, we obtain the growth velocity as a function of temperature $T$, elastic deformation $\bfF_e$, and the state of the system, by identifying the thermomechanical driving force for the melting/freezing process and then assuming a linear kinetic relation between the driving force and growth velocity\footnote{We highlight that the linear kinetic relation is sufficient but not necessary for consistency with the Coleman-Noll procedure.}.
We demonstrate the thermomechanical model and the numerical method by simulating a famous demonstration of regelation: the passage of a weighted wire through a block of ice without splitting it \cite{bottomley1872melting}. 

\subsection{Thermomechanical Model}

We use $\phi=-1$ and $\phi=+1$ to denote the solid ice state and the liquid melt state respectively.
All derivatives, such as gradient and divergence operators, are with respect to spatial coordinates.

\subsubsection{Equilibrium Free Energy}

The Helmholtz free energy is assumed to be composed additively of thermal energy, elastic energy, and interfacial energy. 

\begin{description}

    \item[Thermal energy] 
    
    Building on \cite{penrose1990thermodynamically}, we use the following function for the thermal contribution to the free energy:
    \begin{equation}
        \label{eqn:Free energy-thermal}
        \hat{f}_{Th}(T, \phi) = 
        c_p T \log\left({\frac{T_m^0}{T}}\right) +
        H_l(\phi) L \frac{T_m^0 - T}{T_m^0}
    \end{equation}
    where $L$, $c_p$, and $T_m^0$ are, respectively, the latent heat, heat capacity, and stress-free melting temperature, and are all assumed to be constants; and $H_l(\phi)$ is a smooth version of the sign function such that $H_l(1) = 1$ and $H_l(-1) = 0$.

    \item[Elastic Energy]

    \newcommand{\Ws}{W_\text{solid}}
    \newcommand{\Wm}{W_\text{melt}}

    We use the compressible neo-Hookean form \eqref{eqn: neo-hookean constitutive law} to model the elastic energies of the solid $\Ws$ and the melt $\Wm$, with the shear modulus of the melt being very small but not completely zero for numerical well-posedness.
    Note that the neo-Hookean form in \eqref{eqn: neo-hookean constitutive law} is minimized when the argument is the identity tensor $\bfI$.
    To account for the change in density due to the phase transition in terms of the elastic deformation, we redefine the stress-free state appropriately, and write the elastic free energy density as:
    \begin{equation}
        \label{eqn:Free energy-elasticFinal}
        \hat{f}_{El}(\bfF_e, \phi) = (1-H_l(\phi)) W_\text{solid}
           \left( \left( \frac{\rho^\text{solid}}{\rho^\text{melt}}\right)^{1/d}  \bfF_e \right) + 
           H_l(\phi) W_\text{melt} \left(   \bfF_e \right)      
    \end{equation} 
    where $d$ is the space dimension that we are working in.
   
    To understand \eqref{eqn:Free energy-elasticFinal}, we notice that when $\phi=+1$ corresponding to the liquid state, the first term on the right vanishes, and the argument of $\Wm$ evaluates to $\bfI$ when $\bfF_e = \bfI$.
    That is, in the liquid state, the minimum energy state is when the deformation is the identity (or any rotation); further, the shear modulus is very small, so essentially any volume-preserving deformation is a state of minimum energy.
    When $\phi=-1$, corresponding to the solid state, the second term on the right vanishes, and the argument of $\Ws$ evaluates to $\bfI$ when $\bfF_e = \left( \frac{\rho^\text{solid}}{\rho^\text{melt}}\right)^{-1/d} \bfI$.
    Hence, an isotropic volume change is the state of minimum energy.

    We highlight here that we could consider the process of phase transition as a change in the stress-free deformation $\bfF_{relax}$.
    However, this would lead to a more complex evolution law for $\bfF_e$ as well as require us to formulate the evolution of $\bfF_{relax}$.
    The strategy above has essentially enabled us to use a constant $\bfF_{relax}$ and pose the phase transition process in terms of a change in the elastic energy.

    \item[Interfacial energy] 
    
    As in Section \ref{sec:reg}, to ensure the regularity of the interface for numerical purposes, we consider the following interfacial energy, consistent with the phase field model described in \cite{penrose1990thermodynamically}:
    \begin{equation}
        \label{eqn:Interfacial energy}
        \hat{f}_{In}(\nabla \phi, \phi) =  
        \kappa_1 |\nabla \phi|^2 + \kappa_2 (\phi^2 - 1) ^ 2
    \end{equation}
    where $\kappa_1$ and $\kappa_2$ are small positive modeling constants. 
    
\end{description}

The Helmholtz free energy is then:
\begin{equation}
\begin{split}
    \label{eqn:Helmholtzf}
    \hat f(T,  \bfF_e, \nabla \phi, \phi)  
    &
    = \hat{f}_{Th}(T, \phi) +  \hat{f}_{El}(\bfF_e, \phi)  + \hat{f}_{In}(\nabla \phi, \phi)
    \\
    &
    = c_p T \log\left({\frac{T_m^0}{T}}\right) +
    H_l(\phi) L \frac{T_m^0 - T}{T_m^0}  + 
    (1-H_l(\phi)) W_\text{solid}
    \left( \left( \frac{\rho^\text{solid}}{\rho^\text{melt}}\right)^{1/d}  \bfF_e \right) 
    \\
    & \qquad 
    + H_l(\phi) W_\text{melt} \left( \bfF_e \right) +
    \kappa_1 |\nabla \phi|^2 + \kappa_2 (\phi^2 - 1) ^ 2
\end{split}
\end{equation}

\subsubsection{Kinetics of Melting and Freezing}
\label{sec:growth-velocity}

The evolution of $\phi$ for a given material point due to phase transition occurs by the motion of the ice-water interface, and is modeled using:
\begin{equation}
    \label{eqn: phi_dot pressure melting}
    \dot \phi + \bfv_g \cdot \nabla \phi = 0
\end{equation}
where $\bfv_g$ is the interface velocity.
Note that the material time derivative in Eulerian form is $\dot{\phi} := \parderiv{\phi}{t} + \bfv\cdot\nabla\phi$, and hence \eqref{eqn: phi_dot pressure melting} is identical to \eqref{eqn:phaseFieldFullSeparated}.
Further, with the energy formulated as in \eqref{eqn:Helmholtzf}, the right side of \eqref{eqn: phi_dot pressure melting} being $0$ corresponds to suppressing the nucleation of phases; only existing interfaces are allowed to propagate  \cite{agrawal2015dynamic}.

We now find a form for $\bfv_g$ using the laws of thermomechanics.
Neglecting heat sources, the first and second laws of thermodynamics can be written as:
\begin{align}
    & \rho \dot e = - \divergence \bfq + \bfsigma : \nabla \bfv + \divergence\left(\bft_\phi \dot \phi\right)
    \\
    & \rho \dot s \geq - \frac{\divergence \bfq}{T} + \frac{{}\bfq \cdot \nabla T}{T^2}
\end{align}
where $s$ is entropy per unit mass; $\bfq$ is heat flux; $e$ is internal energy per unit mass.
While most terms above are standard, e.g. \cite{gurtin1982introduction}, we introduce a non-classical work contribution to account for the non-classical phase-field variable.
Specifically, we introduce $t_\phi = \bft_\phi \cdot \bfn$ as the work conjugate of $\phi$ at the boundary, which appears above as $\divergence\left(\bft_\phi \dot \phi\right)$.

Introducing the free energy $f = e - Ts$ and using $f = \hat{f}(T,\bfF_e, \phi, \nabla \phi)$, we follow the classical Coleman-Noll procedure \cite{gurtin1982introduction}, with appropriate extensions for the presence of $\phi$ and $\nabla\phi$, to arrive at the inequality:
\begin{equation}
\label{eqn:entropy-expanded}
    \rho \left(\parderiv{f}{T} + s\right) \dot T + 
    \left(\rho\parderiv{f}{\phi} - \divergence \bft_\phi \right) \dot\phi +
    \left( \rho \parderiv{f}{\nabla \phi}  - \bft_\phi \right) \cdot \nabla \dot \phi+
    \left( \rho \parderiv{f}{\bfF_e} - 
    \rho \left(\nabla \phi \otimes \parderiv{f}{\nabla \phi}\right) \bfF_e^{-T} - \bfsigma \bfF_e^{-T}\right): \dot \bfF_e + 
    \frac{\bfq \cdot \nabla T}{T} \leq 0
\end{equation}

We require each of the five terms in \eqref{eqn:entropy-expanded} to satisfy the inequality independently.
The fifth term gives the classical constraint on $\bfq$, but we assume that $\bfq \equiv \bfzero$ and that the spatial temperature field is constant in time.
Setting the first, third, and fourth terms to $0$ identically gives:
\begin{align}
    \label{eqn:entropy-free-energy}
    s & = -\parderiv{f}{T}
    \\
    \label{eqn:tildeCauchy}
    \bfsigma & = \rho \parderiv{f}{ \bfF_e} \bfF_e^{T}- \rho \left(\nabla \phi \otimes \parderiv{f}{\nabla \phi}\right) 
    \\
    \label{eqn:phi-traction}
    \bft_\phi &= \rho \parderiv{f}{\nabla \phi}
\end{align}
Using \eqref{eqn:entropy-free-energy} and $f = e - Ts$, we obtain the relation $e=c_p T + H_l(\phi) L$ which relates the specific heat and latent heat to the internal energy \cite{penrose1990thermodynamically}; \eqref{eqn:tildeCauchy} is the standard hyperelastic definition of the Cauchy stress tensor \eqref{eqn:hyperelasticMaterial} with an interfacial contribution\footnote{
    The interfacial energy contains spatial derivatives of $\phi$, which leads to a dependence on the deformation gradient and hence the appearance of interface terms in the expression for the stress.
}; and \eqref{eqn:phi-traction} provides an expression for the work-conjugate to $\phi$.

Finally, we are left with the second term in \eqref{eqn:entropy-expanded}:
\begin{equation}
    \left( \rho \parderiv{f}{\phi} - \divergence\left(\rho \parderiv{f}{\nabla \phi}\right) \right) \bfv_g \cdot \nabla \phi \geq 0
\end{equation}
where we have substituted for $\dot \phi$ from \eqref{eqn: phi_dot pressure melting}.
Consequently, we use the following choice for $\bfv_g$, that is sufficient --- but not necessary --- for consistency with thermodynamics:
\begin{equation}
    \label{eqn:vel-driving-force}
    \bfv_g = \kappa \frac{\nabla \phi}{|\nabla\phi|}\left( \rho \parderiv{f}{\phi} - \divergence\left(\rho \parderiv{f}{\nabla \phi}\right) \right)
\end{equation}
where $\kappa > 0$ is a constant.

Evaluating explicitly the variational derivative of $\phi$ --- also denoted as the driving force on $\phi$ --- from \eqref{eqn:vel-driving-force}, by using the expression for $f$ from \eqref{eqn:Helmholtzf}, we have:
\begin{equation}
\begin{split}
    \label{eqn:driving-force-expression}
    & \left( \rho \parderiv{f}{\phi} - \divergence\left(\rho \parderiv{f}{\nabla \phi}\right) \right)    
    \\
    & \qquad = 
    \rho L \frac{T_m^0 - T}{T_m^0} H'_l(\phi) 
    + \rho H'_l(\phi) \left(  
        W_\text{melt} \left(  \bfF_e \right) -  W_\text{solid} \left( \left( \frac{\rho^\text{solid}}{\rho^\text{melt}}\right)^{1/d} \bfF_e \right)
    \right)
    - 2 \kappa_1 \divergence\left( \rho \nabla \phi \right)
    + 4 \rho  \kappa_2 \phi (\phi^2 - 1) 
\end{split}
\end{equation}
We notice that the first 3 terms on the right are localized at the interface, and hence do not drive nucleation, while the last term is decoupled from the deformation and also does not drive nucleation \cite{agrawal2015dynamic}; the model therefore allows phase transitions exclusively by the propagation of existing interfaces.

\subsection{Numerical Simulation}

For the neo-Hookean elastic energy \eqref{eqn: neo-hookean constitutive law}, we use $\mu = 3.52 \mathrm{GPa}, \lambda = 6.54 \mathrm{GPa}, \rho = 900 \mathrm{kg/m}^3$ for the solid corresponding to isotropic ice in the hexagonal phase \cite{greve2009dynamics}; and $\mu = 10^{-3} \mathrm{GPa}, \lambda = 17.8 \mathrm{GPa}, \rho = 1000 \mathrm{kg/m}^3$ for the melt corresponding to liquid water, except that the shear modulus is small but finite for numerical well-posedness.
We use a uniform triangular mesh with element size on the order of $10^{-3}\mathrm{m}$, and use $\kappa_1 = 4 \times 10^{-6}$, $\kappa_2 = 2.3 \times 10^{-4}$, and $\kappa \Delta t = 0.11$ in the appropriate units.

\begin{figure}[htb!]
    \centering
    \includegraphics[width=0.95\textwidth]{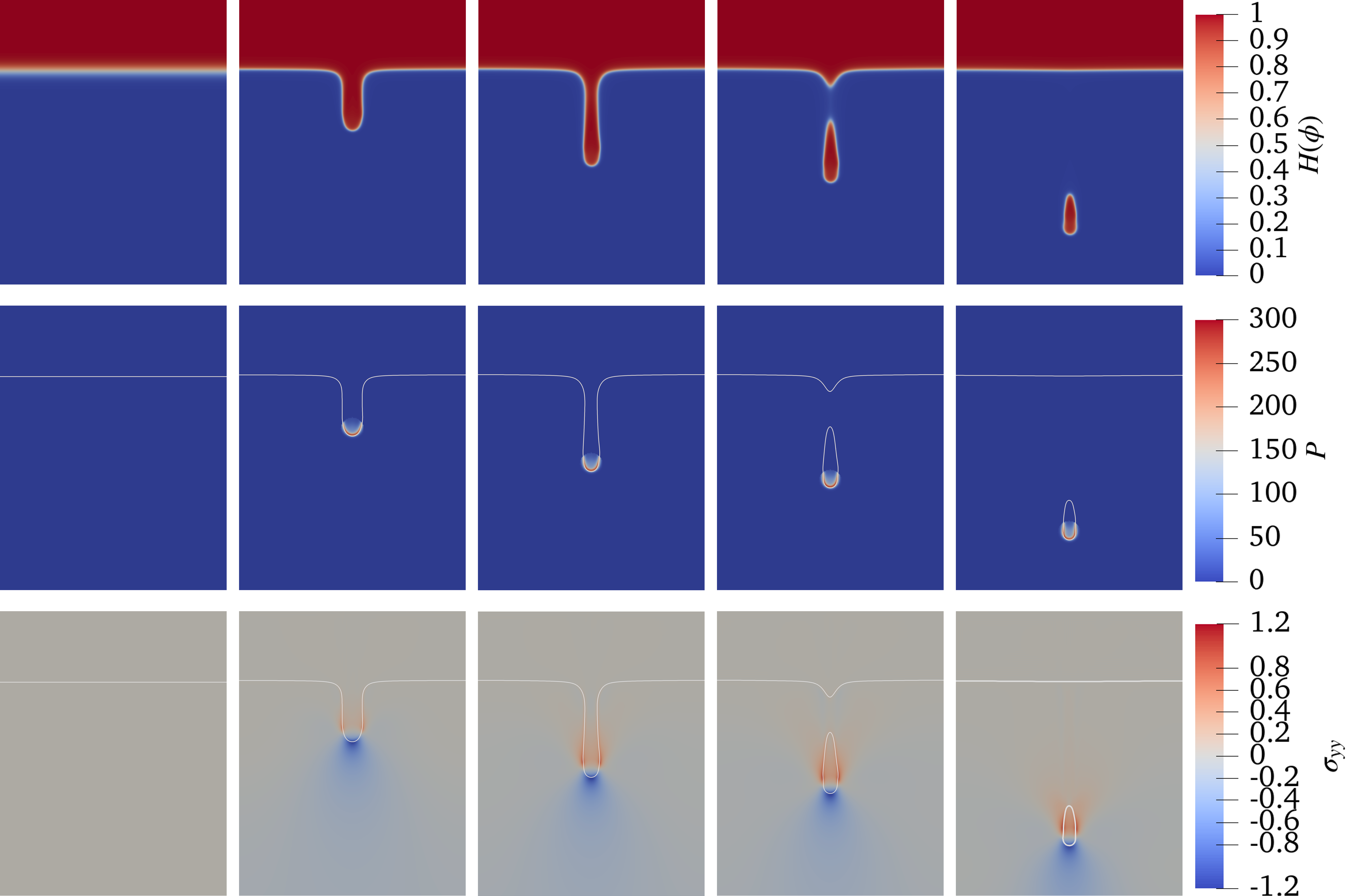}
    \caption{
        Stress-induced melting and refreezing of a solid under a moving load. 
        The top row shows the evolution of $H_l(\phi)$. 
        The middle row shows the load $P$ from \eqref{eqn: P-pressuremelting}, and the white line shows the ice/melt interface, i.e., the levelset  $H(\phi)=1/2$. 
        The bottom row shows the evolution of the vertical normal component of the stress, with the white line showing the interface.  
        From left to right, the columns represent time-steps $0, 35, 65, 80, 130$ respectively.}
    \label{fig:pressure-melting}
\end{figure}

The computational domain is a rectangle $[-0.2, 0.2] \times [-0.25, 0.25] \mathrm{m}^2$.
Initially, the lower portion of the domain is solid ice, and the upper portion is liquid water; Fig. \ref{fig:pressure-melting} (top row) shows the distribution of $\phi$ over time.
Both the ice and water are stress-free at the initial time, with $\bfF_e \approx (1-H_l(\phi)) \sqrt{\frac{\rho^\text{melt}}{\rho^\text{solid}}} \bfI+ H_l(\phi) \bfI$.
The temperature in the water is set to $T_m^0$ implying $L\frac{T_m^0 - T}{T_m^0} = 0$, and in the ice is set below freezing such that $L\frac{T_m^0 - T}{T_m^0} = 0.02$. 
For simplicity, we neglect heat conduction and assume that the spatial temperature does not evolve in time.

We apply a body force $\bfb = -P \bfe_2$ to mimic the traction applied by a weighted wire that presses down on the ice, using
\begin{equation}
\label{eqn: P-pressuremelting}
    P =
    \begin{cases}
         48 H_l'(\phi) \left(1 + \exp\left(-10^4 \left(x^2 + (y-y_0)^2) \right) \right) \right) & \text{ if } \left( x^2 + (y-y_0)^2 \right) < 0.02^2
         \\
         0 & \text{ else }
    \end{cases} 
\end{equation}
In the expression above, $0.02\mathrm{m}$ mimics the radius of the wire, and we set $y_0$ to follow the current location of the wire.
Specifically, we set $y_0$ to be lowest value of the vertical coordinate of the ice/melt interface at the beginning of each time-step.
Further, $H_l'(\phi)$ is only non-zero around the interface, and hence this body force approximates a traction force applied at the interface.

We prescribe the displacement to be $\bfzero$ at the bottom face; prescribe the horizontal displacement to be $0$ and the vertical traction to be $0$ on the left and right faces; and prescribe the traction to be $\bfzero$ at the top face.
When performing the transport update (Section \ref{sec:Advection}), the situations in which we require information from outside $R$ occur only at the top face and  we use the value of the nearest grid point at the boundary.

We note some aspects of the numerical method specific to this problem:
\begin{itemize}
    \item 

    There is no growth on any of the boundaries on which the displacement is prescribed, and hence the growth is unconstrained.
    To smoothly extend the growth velocity over the computational domain $R$, we simply apply the expression in \eqref{eqn:vel-driving-force} over $R$ since the quantities in that expression are defined everywhere in $R$.

    \item 
    From \eqref{eqn: phi_dot pressure melting}, \eqref{eqn:vel-driving-force}, \eqref{eqn:driving-force-expression}, we have that the evolution of $\phi$ is governed by the equation:
    \begin{equation}
    \label{eqn:final-phase-field-evolution}
        \dot{\phi} + \kappa \left( \rho \parderiv{f}{\phi} - \divergence\left(\rho \parderiv{f}{\nabla \phi}\right) \right) |\nabla\phi| = 0
    \end{equation}
    The driving force contains higher-order derivatives that cannot be handled by using a weak form as in standard FEM.
    Consequently, we used a mixed FEM approach. 
    We first evaluate $M:=|\nabla \phi|$ at the beginning of each time-step, and then use an implicit Euler time discretization scheme to solve:
    \begin{equation}
        \dot{\phi} + \kappa \left( \rho \parderiv{f}{\phi} - \divergence\left(\rho \parderiv{f}{\nabla \phi}\right) \right) M = 0
    \end{equation}

    \item 
    Since \eqref{eqn:final-phase-field-evolution} includes higher derivatives that ensure the smoothness of the interface, we do not apply the filtering described in Section \ref{sec:reg}.

    \item 
    We specify the density and elastic deformation of the accreted material --- i.e., the refreezing ice --- assuming that these are identical to the properties of melt immediately before accretion.
    
\end{itemize}

The phase function $H(\phi)$, traction $P$, and the vertical normal component of the stress at different instances of time are shown in Figure \ref{fig:pressure-melting}.
As we can observe from the evolution of phase function, the solid melts in the vicinity of the applied load and refreezes above it, thereby creating a downward-moving melt region within the solid.

\section{Concluding Remarks}
\label{sec:conclusion}

We develop a numerical method to model surface growth using an Eulerian description, following the formulation developed in \cite{naghibzadeh2021surface,naghibzadeh2022accretion}.
This formulation allows the modeling of surface growth without explicitly referring to an unknown time-dependent reference configuration.
We use the phase-field method to enable the use of a fixed discretization within a fixed computational domain for solving the equations numerically. 
    In contrast to Lagrangian descriptions of surface growth, the method proposed in this paper is expected to provide several potential advantages for generality, efficiency, and robustness; however, we highlight that comparing the efficiency and robustness systematically in comparison to other methods is an important future goal.

First, the method in this paper does not impose restrictions on the form of the growth velocity and the stress state of the added material. 
This enables the direct application of the phenomenon of non-normal growth, observed in biological tissues \cite{skalak1997kinematics}, 
Second, it potentially enables the study of the important problem of growth-induced instabilities due to non-uniform growth from a fixed surface, e.g. \cite{abeyaratne2023stability, abeyaratne2022surface, cicconofri2024active}.
Third, using a fixed discretization potentially provides a robust alternative to other available methods which require a moving mesh, re-meshing, and/or extending the computational domain as the computation proceeds.
Fourth, it enables the simulation of the transformation between fluid and solid phases, which is challenging for purely Lagrangian methods that are unsuited for modeling fluids, as well as for Eulerian/Lagrangian hybrid methods that cannot handle the transformation between fluid and solid.

We first applied the method to study non-normal growth observed in some hard biological tissues.
While non-normal growth is not possible to model with typical interface-based methods that use only the normal interfacial velocity, we are able to do this by prescribing the elastic deformation of the accreted material.
Specifically, if the accreted material is pre-stressed at the instant of accretion, it will lead to effective non-normal growth as the pre-stresses relaxes at mechanical equilibrium.
In our work, as a first step we considered only  uniform growth at a flat boundary. 
However, non-uniform growth generates residual stresses in the system that can lead to the development of instabilities at the top surface \cite{abeyaratne2022surface, von2021morphogenesis, abeyaratne2023stability}, and are an important topic for further work.

We next developed a thermomechanical model for first-order phase transformations to describe the freezing and melting of ice under mechanical and thermal loads.
The model does not directly specify the Clausius-Clapeyron relation; instead, this relation emerges as an outcome of a thermomechanical free energy formulation that accounts for general stress states.
While we apply it to the stress-induced melting and refreezing of ice, the overall structure of the model is readily adaptable to various other settings, including ferroelectric \cite{bucsek2019direct} and ferromagnetic \cite{bucsek2020energy} phase transformations, where the extension of the Clausius-Clapeyron relation to general stress states is important \cite{james1986displacive, james1983relation,abeyaratne2006evolution}.
The problem of regelation also involves both fluid and solid phases with transformations in both directions; the deformation of the fluid phase would, in general settings, be exceptionally challenging for a Lagrangian formulation.

In addition to the setting of surface growth, the proposed numerical method can be used to study the mechanics of non-growing solid bodies undergoing extreme deformations where Lagrangian methods perform poorly, e.g. \cite{benson1992computational}.
Specifically, instead of using a moving mesh as in typical Lagrangian and Lagrangian/Eulerian approaches \cite{belytschko2014nonlinear}, transporting the material properties on a fixed discretization following an Eulerian approach is expected to improve the numerical performance by avoiding mesh entanglements under large deformations.
It would also be of interest to apply the Eulerian approach to multiphysics problems, such as poromechanics \cite{clayton2024universal,chua2024deformation,karimi2022energetic} and fracture \cite{clayton2014geometrically,clayton2023phase,agrawal2017dependence,hakimzadeh2022phase, witt2024iga}.

Our approach to solving the transport equations that govern the evolution of the density, elastic deformation, and phase-field used a relatively cumbersome procedure that required the identification of the finite element that contains the transported material point.
While very easy to implement in FEniCS, it requires working directly with nodal values as well as special procedures to handle points that lie outside the computational domain.
An alternative for future work is the Discontinuous Galerkin FEM, that provides an elegant and promising approach that is potentially more efficient \cite{hesthaven2007nodal}.
Similarly, we use a standard FEM formulation that treats the displacement as the primary unknown to solve the incremental equations of hyperelasticity for mechanical equilibrium, but the high elastic contrasts and large deformations could potentially be more accurately and efficiently evaluated by sophisticated mixed approaches following the promising studies in \cite{kadapa2024mixed}.

\section*{}

\paragraph*{Software and Data Availability.}

The code developed for this work and the associated data is available at \\ 
\url{github.com/Kiana-Naghibzadeh/SurfaceGrowth.git}

\paragraph*{Competing Interest Statement.}

The authors have no competing interests to declare.

\paragraph*{Acknowledgments.}

We acknowledge financial support from ARO (MURI W911NF-24-2-0184), NSF (DMREF 2118945, DMS 2108784, DMS 2012259, DMS 2342349), and the Bushnell Fellowship; 
NSF for XSEDE computing resources provided by Pittsburgh Supercomputing Center; 
and Rohan Abeyaratne, Timothy Breitzman, Tal Cohen, and David Rounce for useful discussions.


\newcommand{\etalchar}[1]{$^{#1}$}

\end{document}